\documentclass{article}

\usepackage{PRIMEarxiv}

\usepackage[utf8]{inputenc} 
\usepackage[T1]{fontenc}    
\usepackage{hyperref}       
\usepackage{url}            
\usepackage{booktabs}       
\usepackage{bm}
\usepackage{amsfonts}       
\usepackage{nicefrac}       
\usepackage{microtype}      
\usepackage{lipsum}
\usepackage{fancyhdr}       
\usepackage{graphicx}       
\graphicspath{{media/}}     
\usepackage{epstopdf, booktabs}
\usepackage{subfigure}
\usepackage{diagbox}
\usepackage{caption}
\usepackage[svgnames]{xcolor}
\usepackage{amssymb}

\usepackage{makecell}
\usepackage{pythonhighlight}
\usepackage{listings}
\usepackage[running]{lineno}
\usepackage{setspace}
\usepackage{color}
\usepackage{enumerate}
\usepackage{float}
\usepackage{mathtools, multirow}
\usepackage{threeparttable}
\usepackage{CJKutf8}
\usepackage{lmodern}
\usepackage{amssymb,amsmath}
\usepackage{ifxetex,ifluatex}
\ifnum 0\ifxetex 1\fi\ifluatex 1\fi=0 
  \usepackage[T1]{fontenc}
  \usepackage[utf8]{inputenc}
\else 
  \ifxetex
    \usepackage{mathspec}
  \else
    \usepackage{fontspec}
  \fi
  \defaultfontfeatures{Ligatures=TeX,Scale=MatchLowercase}
\fi

\hypersetup{
colorlinks=true,
linkcolor=blue,
anchorcolor=blue,
citecolor=blue,
filecolor=blue,      
urlcolor=blue,}

\PassOptionsToPackage{usenames,dvipsnames}{color} 
\hypersetup{
            pdftitle={Title},
            pdfauthor={Your name},
            colorlinks=true,
            linkcolor=blue,
            citecolor=Blue,
            anchorcolor=blue,
            filecolor=blue
            urlcolor=Blue,
            breaklinks=true}
\usepackage[authoryear,longnamesfirst]{natbib}
\newcommand{\abs}[1]{\left\lvert #1 \right\rvert}
\newcommand{\pk}[1]{\mathbb{P} \left\{ #1 \right\} }
\newcommand{\expon}[1]{\exp\left(#1\right)}
\newcommand{\QED}{\hfill $\Box$}
\newcommand{\COM}[1]{}

\newcommand{\sign}{\mathrm{sign}}

\newcommand{\mytoprule}{\specialrule{0.1em}{0em}{0em}}
\newcommand{\mybottomrule}{\specialrule{0.1em}{0em}{0em}}

\def\EE{\mathbb{E}}
\def\e{\mathrm{e}} 
\def\I#1{\mathbb{I}\left (#1 \right)}
\def\fracl#1#2{\biggr(\frac{#1}{#2} \biggl) }

\def\IF{\infty}

\def\toweak{\overset{w}\rightarrow}

\usepackage{longtable,booktabs}
\usepackage{graphicx,grffile}
\makeatletter
\def\maxwidth{\ifdim\Gin@nat@width>\linewidth\linewidth\else\Gin@nat@width\fi}
\def\maxheight{\ifdim\Gin@nat@height>\textheight\textheight\else\Gin@nat@height\fi}
\makeatother
\setkeys{Gin}{width=\maxwidth,height=\maxheight,keepaspectratio}
\IfFileExists{parskip.sty}{%
\usepackage{parskip}
}{
\setlength{\parindent}{0pt}
\setlength{\parskip}{6pt plus 2pt minus 1pt}
}
\ifx\paragraph\undefined\else
\let\oldparagraph\paragraph
\renewcommand{\paragraph}[1]{\oldparagraph{#1}\mbox{}}
\fi
\ifx\subparagraph\undefined\else
\let\oldsubparagraph\subparagraph
\renewcommand{\subparagraph}[1]{\oldsubparagraph{#1}\mbox{}}
\fi

\definecolor{c20}{rgb}{0.5,0.5,1}
\definecolor{c30}{rgb}{0.,0.,1.}
\definecolor{c50}{RGB}{34,139,34}
\definecolor{c40}{rgb}{1,0.,0.5}
\definecolor{c60}{rgb}{1,0.9,0.1}
\definecolor{c70}{rgb}{0.50,1.00,0.00}
\def\cl#1{\textcolor{c20}{#1}}
\def\cL#1{\textcolor{c30}{#1}}
\def\kai#1{\textcolor{c40}{#1}}
\def\hu#1{\textcolor{c50}{#1}}

\def\cl#1{#1}
\def\cL#1{#1}
\def\kai#1{#1}
\def\hu#1{#1}

 \numberwithin{equation}{section}
\RequirePackage{colortbl} 
\RequirePackage{soul} 

\newcommand{\tabincell}[2]{\begin{tabular}{@{}#1@{}}#2\end{tabular}}
\newtheorem{theorem}{\bf Theorem\,}[section]
\newtheorem{lemma}{\bf Lemma\,}[section]
\newtheorem{remark}{\bf Remark\,}[section]
\newtheorem{example}{\bf Example\,}[section]

\newtheorem{corollary}{\bf Corollary\,}[section]

\newtheorem{assumption}{Assumption}
\title{Extreme Limit Theory of Competing Risks under Power Normalization
}

\author{
  Kaihao Hu$^{1}$, Kai Wang$^{1}$, 
  Corina Constantinescu$^{2}$, 
  Zhengjun Zhang$^{34}$,  Chengxiu Ling$^{1}$\thanks{Corresponding author: chengxiu.ling@xjtlu.edu.cn}
  \\
  $^{1}$Wisdom Lake Academy of Pharmacy, 
  Xi'an Jiaotong-Liverpool University, SIP, 215123, Suzhou, China 
  \\
    $^{2}$Department of Mathematical Science,  University of Liverpool, L693BX, Liverpool, UK
  \\
  $^{3}$School of Economics and Management, University of Chinese Academy of Sciences, Huairou District,\\ 101418, Beijing, China \\
  $^{4}$Department of Statistics, University of Wisconsin–Madison, 53706, Madison, U.S.A. \\
}

\begin{document}
\maketitle
\begin{spacing}{1.5}
\begin{abstract}
Advanced science and technology provide a wealth of big data from different sources for extreme value analysis.
Classical extreme value theory was extended to obtain an accelerated max-stable distribution family for modelling competing risk-based extreme data in Cao and Zhang (2021).  In this paper, we establish probability models for power normalized maxima and minima from competing risks. The limit distributions consist of an extensional new accelerated max-stable and min-stable distribution family (termed as the accelerated $p$-max/$p$-min stable distribution), and its left-truncated version. The consistency and asymptotic normality are obtained for the maximum likelihood estimation of the parameters involved in the accelerated $p$-max and $p$-min stable distributions when it exists. The limit types of distributions are determined principally by the sample generating process and the interplay among the competing risks, which are illustrated by common examples.  Further, the statistical inference concerning the maximum likelihood estimation and model diagnosis of this model was investigated.  
Numerical studies show first the efficient approximation of all limit scenarios as well as its comparable convergence rate in contrast with those under linear normalization, and then present the maximum likelihood estimation and diagnosis of accelerated $p$-max/$p$-min stable models for simulated data sets. 
Finally, two real datasets concerning annual maximum of ground level ozone and survival times of Stanford heart plant demonstrate the performance of  our accelerated $p$-max and accelerated $p$-min stable models.
\end{abstract}

\end{spacing}

\newpage
\section{Introduction}

In classical extreme value theory, the central result is the Fisher-Tippett-Gnedenko theorem which specified the form of the limit distribution for centered and normalized maxima \citep{leadbetter2012extremes}. Let $X_{1}, X_{2}, \dots, X_{n}, \ldots$ be a sequence of independent and identically distributed (i.i.d.) random variables with common distribution function (d.f.) $F$ 
and denote by $M_n = \max (X_{1}, X_{2}, \dots, X_{n})$ the sample maxima.
If there exist some normalization constants $a_{n}>0, b_{n}\in\mathbb R$ and a non-degenerate d.f. $G$ such that
\begin{equation}\label{l-MDA}
    \mathbb{P}\left ( a_{n}\left ( M_{n}-b_{n} \right )\le x  \right ) \toweak G\left ( x \right )\quad \mbox{as}\ n\to\IF,
\end{equation}
then $G$ must be {one of the following three} $l$-types (i.e., $G(ax+b)$ is of the same $l$-type of $G(x)$ for any $a>0, b\in\mathbb R$):   
\begin{eqnarray*}
   && \mbox{Gumbel}: \Lambda(x)=\expon{-\e^{-x}  },\quad x\in\mathbb{R};\\
   &&\mbox{Fr\'echet}: \Phi_\alpha(x)= \begin{cases}0, & x \leq 0 , \\ \expon{-x^{-\alpha} }, & x>0;\end{cases} \\
  &&\mbox{Weibull}: \Psi_\alpha(x)= \begin{cases}\expon{ -\left ( -x \right )^{\alpha }} , & x \leq 0, \\ 1, & x>0. \end{cases}
\end{eqnarray*}
Here $\alpha$ is a positive parameter and $\toweak$ stands for the weak convergence (convergence in distribution). We denote by $F\in D_l(G)$ if Eq.\eqref{l-MDA} holds, indicating that $F$ belongs to the max-domain attraction of  $G$  under linear normalization.
Note that the class of limit distributions $G$ coincides in the so-called the max-stable distribution under linear normalization or simply \textit{$l$-max stable} d.f., since for any $n = 2, 3, \ldots$, there are some constant sequences  $a_n > 0$ and $ b_{n}\in\mathbb R$ 
such that $G^n(a_n^{-1}x+b_n) = G(x)$. 
In addition, a unified form of the three types of $l$-max stable distributions is the \textit{generalized extreme value distribution} (GEV) \kai{$G_\xi(x; \mu, \sigma)$} given by 
\begin{eqnarray}\label{GEV}
    G_\xi(x; \mu, \sigma) = \expon{-\left[1+\xi\left(\frac{x-\mu}{\sigma}\right)\right]_+^{-1/\xi}},
\end{eqnarray}
where $\xi\in\mathbb R$, is the so-called extreme value index, measuring the tail heaviness of $F$. The $G_\xi$ with $\xi >, =, <0$ corresponds to respectively the Fr\'echet, Gumbel and Weibull 
distributions, and $\alpha = 1/|\xi|$. Here the case of $\xi=0$ can be understood as $\lim_{\xi \to 0} G_{\xi}$. 
\\
As shown by \citet{Mohan1993}, the class of distributions $F\in D_l(G)$ can be essentially extended to a wider class of $F\in D_p(H)$  with $H$ a \textit{$p$-max stable} distribution \citep{pantcheva1984}, i.e., there exist some normalization constants $\alpha_{n}, 
\beta_{n}>0$ and a non-degenerate d.f. $H$ such that 
\begin{equation}\label{p-MDA}
   \mathbb{P}\left \{ \alpha_{n} \left | M_ {n} \right |^{\beta _{n}}\sign\left ( M_{n} \right )  \le x  \right \} \toweak H(x),
\end{equation}
with $\sign(x)=-1,0,1$ for $x<,=,>0$. Then $H$ must be one of the following six $p$-types (i.e., $H(A|x|^{B}\sign(x))$ is of the same $p$-type of $H(x)$ for any $A,B>0$)
\begin{eqnarray*}
    H_{1,\alpha}\left ( x \right )&=&\begin{cases}
    0,& x\le 1{,}\\ \Phi_\alpha(\log x),& x>1;
    \end{cases}\\
    H_{2,\alpha}\left ( x \right ) &=&\begin{cases}
    0,&x\le0,\\ \Psi_\alpha(\log x),& 0<x<1,\\1,&x\ge 1
    ;
    \end{cases}\\
 H_{3,\alpha}\left ( x \right) &=& \begin{cases}
    0{,} &x\le -1{,}\\
     \Phi_\alpha(-\log(-x)),&-1<x< 0, \\ 1,&x\ge 0;
          \end{cases}     \\
     H_{4,\alpha}\left ( x \right )&=&\begin{cases}
    \Psi_\alpha(-\log(-x)) ,&x< -1, \\ 1,&x\ge -1;
    \end{cases}\\
     H_{5}(x) &= & \Phi_1(x);\\
    H_{6}(x) &=& \Psi_1(x).
    \end{eqnarray*}
Here $\alpha$ is a positive parameter. Note that an alternative form of  a  sample maxima $M_n$ from risk $X\sim F$, is given by 
\begin{eqnarray}\label{sample-min-max}
M_n =-\min_{1\le i \le n}(-X_i) \stackrel{d}{=} - \underline M_n,    
\end{eqnarray}
where $\stackrel d=$ stands for equality in distribution, and the sample minima $\underline M_n$ is generated from $-X\sim \underline F(x) := 1-F(-x-0)$. Thus, the dual limit theorems of sample minima from risk $\underline F$ under linear and power normalization  hold with non-degenerated  distributions $\underline G(x)$ and
$\underline H(x)$, if and only if $F\in D_l(G)$ and  $F\in D_p(H)$, respectively. We call $\underline{G}$ and $\underline H$ as \textit{$l$-min stable} and \textit{{$p$}-min stable} distributions, see \citet[Theorem 1.8.3]{leadbetter2012extremes} and \citet[Definition 2]{Grigelionis2004}.

It is well known  that the selection of normalization constants $(a_n, b_n)$ (or $(\alpha_n, \beta_n)$) are completely determined by the tail behavior of $F$, see e.g., \citet[Theorems 2.1--2.6 and Theorem 3.1]{Mohan1993} and \citet[Corollary 1.6.3]{leadbetter2012extremes}. 
Each pair of these normalization constants can be linearly (or powerly) equivalent, i.e.,  
\begin{eqnarray*}
   &&a_n' \approx a_n /a, \quad b_n' \approx b_n + b/a_n\quad \mbox{for}\ a>0,\, b\in\mathbb {R}; \\
   &&\alpha_n' \approx (\alpha_n/A)^{1/B}, \quad \beta_n' \approx \beta_n/ B\quad \mbox{for}\ A, \, B>0.
\end{eqnarray*}
These equivalent normalization constants lead to the same type of non-degenerate limit distributions in Eqs.\eqref{l-MDA} and \eqref{p-MDA}, i.e., $G(ax+b)$ and $H\big(A|x|^{B}\sign(x)\big)$, respectively. 
We refer to \citet[Theorem 1.2.3 (Khintchine)]{leadbetter2012extremes} and \citet[Lemma 2.1]{Baracat2002} for detailed discussions.

There are many studies on the two normalized extremes and its relationships as well as its applications. \citet{Mohan1993} obtained equivalent conditions of two max-domain attractions (MDAs) and its relationship between the linear and power normalization constants. The convergence rate relationship of extremes under linear and power normalization was studied further by \citet{Barakat2010}, see e.g.,  
\citet{Chen2012, GengLi2016} and \citet{LIAO201440} for related 
studies of power normalized extremes from specific risk $F$, following general error distribution, logarithm general error distribution and  skew-normal  distribution, respectively.  Meanwhile, asymptotic behavior  of extremes such as minima and extreme order statistics under power normalization was examined by \citet{Grigelionis2004} and \citet{PengZuoxiang2012}, respectively. In addition, \citet{Silvestrov1998} and \citet{Baracat2002} considered random sample size and investigated asymptotic behavior of extremes under linear and power normalization, see e.g., 
 \citet{Mavitha2016, Barakat2020} for recent analysis of exponential normalized extremes, with flexibility in the sign of extremes compared with power normalization. It has also been extended into the extreme order statistics, see \citet{barakat2023limit}.

One common assumption in the aforementioned studies is the identical distribution of the data, which may restrict its applications since extreme data arising in most real-world problems may cluster in space and across time, e.g., ground level ozone data \citep{Fuentes2003}. To handle big data with complex structures, such as competing risks in finance and economics fields, \citet{cao2021} established the new extreme value theory of maxima of maxima under linear normalization, which helps decomposing systemic risk into competing risks for learning risk patterns and better risk management \citep{Ji2021}. Specifically, assume that $X_1, X_2, \ldots, X_n$ is a collected dataset from $k$ independent sources (units, regions, periods etc) $F_j, j=1,\ldots, k$, from which a random sample  $X_{j,i}, i=1, \ldots, n_j$ comprise the mixed sample $X_i, 1\le i\le n$ with $n=n_1+n_2+\dots+n_k$. Denote by $M_{j,n_j} = \max_{1\le i\le n_j} X_{j,i}$  the sample maximum of the $j$-th source with sample size $n_j$. Then the sample maxima $M_n$ can be rewritten as 
\begin{eqnarray} \label{MaxMax}
    M_n = \max_{1\le j \le k} M_{j, n_j}.
\end{eqnarray}
Note that the extreme behavior of $M_n$ drawn in a competing risk scenario is completely different from that of maxima in multivariate setting, namely, the  multivariate extreme value theory \citep{resnick2008extreme}. We will see that the limit behavior will depend on both the sample frequency and the risk heaviness. 
The limit distribution of the linear normalized extreme of Eq.\eqref{MaxMax} is given by the so-called  accelerated $l$-max stable distributions under mild conditions on the sample length and data generating sources \citep[Theorem 2.1]{cao2021}, i.e., $G(x)=\prod_{j=1}^{k'}G^{(j)} (x)$ for some positive integer $k'$ where $G^{(j)}$ is a $l$-max stable d.f.. 

Given the advantages of power normalization over linear normalization in terms of the convergence rate and the wider applications  of limit distributions \citep{Barakat2010}, and the common data structure caused by competing risks, this paper aims to answer the following questions.
\begin{itemize}
    \item[A.] What is the limit distribution of heterogeneous extremes under power normalization?
    \item[B.] What are the similarities and differences of the behavior of heterogeneous extremes under linear and power normalization?
\end{itemize}

The main focus of the paper is to develop the extreme type theory of heterogeneous extremes under power normalization. The new limit distributions consist of the so-called accelerated $p$-max stable distributions that may be the $p$-max stable distributions multiplied (accelerated by) other competing $p$-max stable distributions or its product, 
and its left-truncated version 
introduced in Theorem \ref{TheoremMaxMax}, extending both $p$-max stable and accelerated $l$-max stable family. The limit distributions of minima of competing risks are established as well as the statistical inference concerning maximum likelihood estimation (mle) of this new family are discussed, see Sections \ref{minima} and \ref{sec:mle}. Specifically, we obtained the consistency and asymptotic normality of the mle in Theorems \ref{Theoremconsistency} and  \ref{Theoremasymptotic normality}.
We give further common examples in Section \ref{sec: Examples} to show how the limit distributions given in Theorem \ref{TheoremMaxMax}, representing the three different scenarios, namely, the accelerated $p$-max stable case, the single $p$-max stable case and the left-truncated case, are determined by the sample generating process and interplay among the competing risks. Roughly speaking, the three cases arise when one of the limit distributions is non-degenerate, degenerate at a point in the bound or interior of the support of another limit distribution.  Numerical analysis conducted in Section \ref{sec: simulation} illustrates first the efficiency of the distribution approximations established by Theorem \ref{TheoremMaxMax}, followed by the comparison of the convergence rate of competing risks under linear and power normalization by employing the Kolmogorov-Smirnov (KS) test, Cram\'{e}r-von Mises (CVM) test and Anderson-Darling (AD) test. These model goodness of test approaches, together with the maximum likelihood estimation mechanism developed in Section \ref{sec:mle}, were further used to detect the competing risks for both simulated data sets and real data sets concerning annual maximum ground-level ozone and survival times of Stanford heart plant patients. These results demonstrate the wide applications of our models. Many profound studies in this direction are worth doing in the future, e.g., the exponential normalized extreme of competing risks, and extreme behaviour of competing risks with random sample size, see e.g., \citet{Baracat2002} and \citet{Barakat2020}. 

The remainder of the paper is organized as follows. Section \ref{sec: MainResults} presents the main results. Discussions of limit theorems of minima of competing risks and the asymptotic properties of maximum likelihood estimation involved are given in Section \ref{sec:discussion}. Section \ref{sec: Examples} provides illustrating examples, followed by numerical studies in Section \ref{sec: simulation}.
Finally, an application of two real data concerning ground-level ozone and survival times of Stanford heart plant patients is demonstrated in Section \ref{sec:Real data}. The proofs and relevant 
lemmas are deferred to Appendix.

\section{Main Results}\label{sec: MainResults}
In what follows, we suppose that all samples are independent unless 
stated otherwise. Further,  for a d.f. $F$, we denote by $\gamma (F) = \sup\{x \in \mathbb R: F(x)<1\}$ for the right-endpoint. 
We investigate in Theorem \ref{TheoremMaxMax} below the limit behavior of competing extremes (see Eq.\eqref{MaxMax}) under the following mild conditions on each sub-block maximum. 

Let $M_{j,n_j} = \max_{i=1}^{n_j} X_{j,i}$ be the sample maxima from d.f. $F_j, \ j=1,  2, \ldots, k$. Assume that there exist some normalization constants $\alpha_{j,n_j}, \ \beta_{j,n_j}>0$ and a non-degenerate distribution $H^{(j)}$ such that
\begin{equation}\label{ComponentPowerLimit}
     \mathbb{P}\left ( \alpha_{j,n_{j} }\left | M_{j,n_j} \right |^{\beta _{j,n_{j} }}\sign (M_{j,n_j}) \le x    \right ) {\toweak} H^{(j)}\left ( x \right )
\end{equation}
as $n_j \to\infty$. That is, each risk $F_j\in D_p(H^{(j)})$, satisfying Eq.\eqref{p-MDA} with the above normalization constants and limit $H^{(j)}$ being 
one of the six $p$-types. We focus on the case with $k=2$, the maxima of two competing risks, and the general case can be analysed similarly. Denote further 
\begin{equation} 
    \label{NormingMaxMax}
\alpha_{n}=\alpha_{1,n_{1}}\left ( \frac{1}{\alpha_{2,n_{2}}}  \right )^{\frac{\beta_{1,n_{1}}}{\beta_{2,n_{2}}}}, \quad \beta_{n}=\frac{\beta_{1,n_{1}}}{\beta_{2,n_{2}}}.
    \end{equation}
We will see below that the limit behavior of $\alpha_n$ and $\beta_n$ plays a key role in the determination of limit types of the competing risks. Similar to \citet[Theorem 2.1]{cao2021}, we consider first mixed and dominated relationships between two risks to derive the limit distributions of maxima of maxima under power normalization. Different from linear normalization, a new left-truncated limit was established.

\begin{theorem}
    \label{TheoremMaxMax}
        Let $M_n = \max(M_{1, n_1}, M_{2, n_2})$ be the max of two sample maxima. Suppose that condition \eqref{ComponentPowerLimit} holds for $M_{j, n_j}$ with normalization constants $\alpha_{j,n_j}, \beta_{j, n_j}, j=1,2$.  Assume that there exist two constants $A, B\in[0,\infty]$ such that 
    \begin{equation}
            \label{ConditionNorming}
            \alpha_n \to A,\quad \beta_n \to B         \end{equation}
as $\min(n_1, n_2) \to\infty$. 

    \begin{itemize}
\item[(i)] If $AB\in(0,\infty)$, i.e., condition \eqref{ConditionNorming} holds with two positive constants $A$ and $B$, then 
$$\pk{ \alpha _{2,n_{2} }\left | M_{n}  \right |^{\beta _{2,n_{2} } }\sign\left ( M_{n} \right )\le x } \toweak H^{(1)}\big( A|x|^{B} \sign(x) \big) H^{(2)}(x).$$
\item[(ii)] 
The following limit distribution holds 
$$
\pk{ \alpha _{2,n_{2} }\left | M_{n}  \right |^{\beta _{2,n_{2} } }\sign\left ( M_{n} \right )\le x } \toweak  H^{(2)}(x)
$$
provided that one of the following four conditions is satisfied.
\begin{itemize}
\item[a).] $H^{(2)}$ is one of the same $p$-types of $H_{1,\alpha},H_{2,\alpha},H_{5}$, and $H^{(1)}$ is one of the same $p$-types of $H_{3,\alpha},H_{4,\alpha},H_{6}$.
\item[b).] $H^{(2)}$ is one of the same $p$-types of $H_{1,\alpha},H_{2,\alpha},H_{5}$, and $H^{(1)}$ is one of the same $p$-types of $H_{1,\alpha},H_{5}$. In addition, Eq.\eqref{ConditionNorming} holds with $A=\infty$ and $0 \le B < \infty$. 
\item[c).] $H^{(2)}$ is one of the same $p$-types of $H_{1,\alpha},H_{2,\alpha},H_{5}$, and $H^{(1)}$ is the same type of $H_{2,\alpha}$. In addition, Eq.\eqref{ConditionNorming} holds with  $\gamma(H^{(1)}) \le A < \infty$ and $B=0$ or $A = \infty$ and $\ 0 \le B < \infty$. 
\item[d).] Both $H^{(1)}$ and $H^{(2)}$ are one of the same $p$-types of $H_{3,\alpha},H_{4,\alpha},H_{6}$. In addition, Eq.\eqref{ConditionNorming} holds with $0\le A\le -\gamma(H^{(1)})$ and $ B=0$ or $A=0$ and $ \ 0\le B <\infty$. 
\end{itemize}
\item[(iii)] The following limit distribution holds 
$$
\pk{ \alpha _{2,n_{2} }\left | M_{n}  \right |^{\beta _{2,n_{2} } }\sign\left ( M_{n} \right )\le x } \to  H^{(2)}(x) \I{x > x_0}
$$
if $\alpha_n$ and $\beta_n$ given in Eq.\eqref{NormingMaxMax} satisfy  $\lim_{n \to \infty} (\log \alpha_n)/\beta_n=-C<0$ (in this case Eq.\eqref{ConditionNorming} holds with $A=0, B=\infty$), 
and both $H^{(1)}$ and $H^{(2)}$ are one of the same type of $H_{1,\alpha},H_{2,\alpha},H_{5}$, or one of  the same type of $H_{3,\alpha},H_{4,\alpha},H_{6}$. The jump point $x_0$ equals $\e^C, -\e^C$ for the two cases, is supposed to be in the interior of the support of $H^{(2)}$.

 \end{itemize}
\end{theorem}

Note that Theorem \ref{TheoremMaxMax}(i) can be directly verified by Lemma \ref{lemmaA} \citep[Lemma 2.1]{Baracat2002}. It shows that the positive limits in Eq.\eqref{ConditionNorming} ensure $M_{1, n_1}$ under the alternative normalization constants $\alpha_{2,n_2}, \beta_{2, n_2}$ to converge in distribution to the same $p$-type max stable distribution. In addition, the limit in Theorem \ref{TheoremMaxMax}(i) is the product of two extreme value distributions, $H^{(1)}\big( A|x|^{B} \sign(x) \big) H^{(2)}(x)$. Although it is of product form, it can still be reduced to the six types of extreme value distributions in some cases. For instance, $H_{1, \alpha}(x^B)H_{1, \alpha}(x) = \expon{-(1+B^\alpha)(\log x)^{-\alpha}},\, x>1$ is again of the $p$-type of $H_{1,\alpha}$. However, in some situations, the limit product form can not be reduced to any one of the six types of extreme value distributions, see e.g., Example \ref{p-Norm-logFrechet}. In general, we call $H$ the accelerated $p$-max stable distribution if it can be written as a product of $p$-max stable distributions $H^{(j)}, \, 1\le j \le k$: 

\begin{eqnarray}\label{accelerated-p}
    H(x)=H^{(1)}(x)H^{(2)}(x)\cdots H^{(k)}(x). 
 \end{eqnarray} 
 
Apparently, the accelerated $p$-max stable distributions family includes the $p$-max stable distributions as special cases, and it is larger 
than the accelerated $l$-max stable distributions given by \citet[Definition 2.1]{cao2021}. 
In real-world applications, the accelerated cases can be regarded as accelerated risks from a competing scenario.
For example, in a medical setting, a patient with several diseases may enter the emergence room due to the deterioration of the most severe disease. Denote by the most severe risk $Z = \max(Z_1,\ldots, Z_k)$ with risk $Z_j \sim H^{(j)}$ associated with disease $j$ independently. We have the accelerated risk $Z$ following the accelerated $p$-max stable distribution $H(x)$ given by Eq.\eqref{accelerated-p}. Clearly, we have

$$
H(x) \le \cl{\min_{1\le j\le k}H_{-j}(x)}, 
$$

where $H_{-j} = \prod_{1\le i\neq j \le k} H^{(i)}$. It indicates that the severity of the patient with $k$ different diseases is higher than other patients with $k-1$ or less diseases under the same risk level probability.
This situation can be also implied in terms of the emergence room entrance probability (ERP) below  

$$ERP = 1-H(x)=1-H^{(1)}(x)H^{(2)}(x)\cdots H^{(k)}(x) \ge \max_{1\le j\le k}(1-H^{(j)}(x)).$$

That is, 
the ERP of a patient with multiple/compound diseases is higher than that of a patient with any individual disease. In survival analysis, the survival time $T_j$ associated with $j$-th disease is the primary quantity for study, thus the survival time of the accelerated risk $Z = \max_{1 \le j \le k} Z_j$ 
is determined by
$T=\min_{1\le j\le k} T_j$, which survival probability of all diseases 
becomes smaller for more severe patients and thus the lifetime is shortened most due to the accelerated risk. The accelerated risk can also be regarded as the systemic risk being accelerated from individual risk given a fixed confidence interval in financial scenarios \citep{cao2021}. Further applications include the endopathic and exopathic dynamic competing risks model, providing contributions to risk measurement in financial markets, environmental disasters and global disease pandemics \citep{Ji2021}.

\begin{remark} \label{remark1}
a) Different from Theorem \ref{TheoremMaxMax}(i) for the accelerated limit distributions, the limit $H^{(2)}$ in Theorem \ref{TheoremMaxMax}(ii) can be roughly regarded as the case that the maxima $M_{2, n_2}$ is larger than $M_{1, n_1}$ in probability. Specifically, 
\begin{itemize}
    \item[(i)] case a) is apparent with $M_{2,n_2}$ greater than $M_{1, n_1}$
since  $X_{2, i}> 0\ge \gamma(F_1)\ge X_{1, i'}$, implying thus $M_n = M_{2, n_2}$ almost surely.
\item[(ii)] the cases of b), c) and d) are obtained when $x_n = \alpha_n|x|^{\beta_n} \to x_\infty\ge \gamma(H^{(1)})$ for any given $x$ in the support of $H^{(2)}$ since elementary calculations indicate that \begin{eqnarray*}
  &&\pk{\alpha_{2, n_2}|M_{1, n_1}|^{\beta_{2, n_2}} \sign(M_{1, n_1}) \le x} \\&= &\pk{\alpha_{1, n_1}|M_{1, n_1}|^{\beta_{1, n_1}} \sign(M_{1, n_1}) \le x_n},  
\end{eqnarray*}
see e.g.,  Examples \ref{FandF}, \ref{WandW}, \ref{UandU}, \ref{ex3.8} and \ref{ex3.9}.
\end{itemize}

b) The left-truncated version of $H^{(2)}$ corresponds to the limit distribution in Theorem \ref{TheoremMaxMax}(iii), which may occur when $x_n = \alpha_n|x|^{\beta_n} \to x_\infty$ with
$H^{(2)}$. For instance,  

$$
x_\infty = \left\{ 
\begin{array}{ll}
\gamma(H^{(1)}), &x>x_0,\\
\gamma_*(H^{(1)}) := \inf\{x\in\mathbb R: H^{(1)}(x)>0\}, & x \le x_0
\end{array}
\right.
$$

with $x_0$ being in the interior of the support of $H^{(2)}$, see e.g.,  Examples \ref{ex33} and \ref{GandF}.
\end{remark}

\section{Discussions}\label{sec:discussion}
We will discuss first the limit theorems for minima of competing risks and then statistical analysis of our models including  maximum likelihood estimations of the parameters involved.

\subsection{Limit Theorems of Minima of Competing Risks}\label{minima}
In many practical applications, of interest is to study the minima of competing risks. For instance, in survival analysis, the survival time $T$ of an patient with multiple diseases is determined as the life time of the patient surviving from each of the diseases. Recalling Eqs.\eqref{sample-min-max} and \eqref{MaxMax}, we have for sample minima $\underline M_{j,n_j}$ from $k$ independent parents $-X_j \sim \underline F_j,\, j=1, \ldots, k$ 
$$\underline M_n = \min_{1\le j \le k} \underline M_{j,n_j} = - \max_{1\le j \le k} M_{j,n_j}$$
with sample maxima $M_{j,n_j}$ from $X_j \sim F_j$. Since $A|(-x)|^B \sign(-x) = - A|x|^B \sign(x)$, we have the same power normalization constants for both $\underline M_{j,n_j}$ and $M_{j,n_j}$ when considering its limit behavior. 
\begin{remark}\label{remark2} 
The risk $X_j$ follows $F_j \in D_p(H^{(j)})$ with sample maxima $M_{j,n_j}$ satisfying Eq.\eqref{p-MDA} with normalization constants $(\alpha_{j,n_j}, \beta_{j,n_j})$ if and only if sample minimum $\underline M_{j,n_j}$ from risk $-X_j \sim \underline F_j$ satisfies 
\begin{equation}\label{ComponentPowerLimitmin}
\begin{aligned}
       \mathbb{P}\left ( \alpha_{j,n_{j} }\left | \underline{M}_{j,n_j} \right |^{\beta _{j,n_{j} }}\sign (\underline{M}_{j,n_j}) \le x    \right ) {\toweak} \underline{H}^{(j)}\left ( x \right ),
\end{aligned}
\end{equation}
where $\underline{H}^{(j)}(x)=1-H^{(j)}(-x-0)$ is the so-called $p$-min stable distribution, see \citet[Corollary 1]{Grigelionis2004}. 
\end{remark}
\begin{corollary}\label{Cor-1}
Assume that condition \eqref{ComponentPowerLimitmin} holds for independent risks $-X_j \sim \underline F_j, j=1,2$, and condition \eqref{ConditionNorming} is satisfied for some constants $A$ and $B$, then all claims in Theorem \ref{TheoremMaxMax} hold for minima of competing risks $\underline M_n$ and limit distribution $\underline H(x) = 1- H(-x-0)$ with $H$ the accelerated $p$-max stable distribution and its left-truncated version given in Theorem \ref{TheoremMaxMax}. 
\end{corollary}
 In parallel with accelerated $p$-max stable distributions, we call $\underline H$ as an \textit{accelerated $p$-min stable distribution}, if $H(x) = 1-\underline H(-x-0)$ is an accelerated $p$-max stable distributions.

\subsection{Maximum Likelihood Estimation} \label{sec:mle}
Following our main results for the limit distribution approximations of extremes from multiple competing sources, two natural questions arising are as below. 
\begin{itemize}
    \item[a).] How to estimate the relevant parameters involved in limit distributions in Theorem \ref{TheoremMaxMax} if we give a sample of extremes from competing risks?
    \item[b).] How to detect if extremes can be modeled by our models, namely, data are from multiple competing risks?  
\end{itemize}
In what follows, we will focus on the statistical analysis by addressing the aforementioned questions subsequently. 

Without loss of generality, we suppose that $M_n = \max(M_{1,n_1}, M_{2,n_2})$ with $M_{1,n_1}, M_{2,n_2}$ and $x$ being of the same sign at infinity. Thus, 
\begin{eqnarray}\label{max-no-norming}
    \pk{M_n \le x} 
    &=&\pk{M_{1,n_1} \le x}\pk{M_{2,n_2} \le x} \notag\\
    &=& \pk{\alpha_{1,n_1} |M_{1,n_1}|^{\beta_{1,n_1}} \sign(M_{1,n_1}) \le\alpha_{1,n_1} |x|^{\beta_{1,n_1}} \sign(x)} \notag\\
    &&\times\pk{\alpha_{2,n_2} |M_{2,n_2}|^{\beta_{2,n_2}} \sign(M_{2,n_2}) \le \alpha_{2,n_2} |x|^{\beta_{2,n_2}}\sign(x)} \notag \\ 
    &=& H^{(1)}(\alpha_{1,n_1} |x|^{\beta_{1,n_1}}\sign(x)) H^{(2)}((\alpha_{2,n_2} |x|^{\beta_{2,n_2}}\sign(x)) \notag\\
    &=:&\widetilde H^{(1)}(x) \widetilde H^{(2)}(x).
\end{eqnarray}
We see that Eq.\eqref{max-no-norming} indicates that the practical applications of our models are free of power normalization constants, as they are absorbed in $\widetilde H^{(1)}(x)$ and $\widetilde H^{(2)}(x)$,  and thus 
they are free of sample sizes,
see also \citet{Nasri1999}. 

In addition, the six types of $p$-max stable distributions can be written into two generalized forms, in accordance with the support being in $(-\infty,0)$ and $(0,\infty)$, respectively. In the following, we focus first on  the non-negative risks we encountered in most application-oriented fields.  Specifically, the distributions,  being of the same type of $H_{1,\alpha}, H_{2,\alpha}$ or $H_{5}$, can be uniformly written as 
$H^{\xi}_{1}(x; \mu, \sigma) = G_\xi(\log x;\mu,\sigma)$ with $G_\xi$ given by  Eq.\eqref{GEV}, see e.g.,  \citet{Nasri1999}.

Let $z_1, \ldots, z_m$ be an i.i.d. observed data set from parent $Z\sim H_1^\xi(z;\mu,\sigma)$. The log-likelihood function is thus given by
$\ell_1(\boldsymbol z;\mu,\sigma,\xi) = \sum_{i=1}^m \log(H_1^\xi(z_i;\mu,\sigma))'$. 

\COM{\hu{In addition, if $Z \sim \underline{H}_2^\xi$, a $p$-min stable distribution function, the log-likelihood function is given by}
\begin{equation}\label{log-likelihood-p-min}
\ell_2(\boldsymbol z;\mu,\sigma,\xi)=
\left\{ 
\begin{array}{ll}
-\sum_{i=1}^{m}\left[1-\xi\fracl{\log z_i-\mu}{\sigma}\right]^{-1/\xi}-\left(1+\frac{1}{\xi}\right)\sum_{i=1}^{m}\log\left[1-\xi\left(\frac{\log z_i-\mu}{\sigma}\right)\right]     &  \\
- m\log \sigma -\sum_{i=1}^{m}\log  z_i, & \xi\neq0;\\
-\sum_{i=1}^{m} \expon{ \frac{\log z_i-\mu}{\sigma}}+\sum_{i=1}^{m} \frac{\log z_i-\mu}{\sigma} - m\log \sigma -\sum_{i=1}^{m}\log z_i. & \xi=0.
\end{array}
\right.
\end{equation}
}
Further, if $Z\sim H$,  
an accelerated $p$-max stable distribution as  
$H(z; \boldsymbol\theta) = H_1^{\xi_1}(z;\mu_1,\sigma_1)H_1^{\xi_2}(z;\mu_2,\sigma_2),$  
then the log-likelihood function $\ell(\boldsymbol\theta)$ with $\boldsymbol\theta = (\mu_1, \sigma_1,\xi_1, \mu_2, \sigma_2,\xi_2)$ is given by
\begin{eqnarray}\label{log-likelihood-acc-pmax}
\ell(\boldsymbol\theta) = \sum_{i=1}^m \log\left[H_1^{\xi_2}(z_i; \mu_2, \sigma_2)\e^{\ell_1(z_i;\mu_1,\sigma_1, \xi_1)} + H_1^{\xi_1}(z_i; \mu_1, \sigma_1)\e^{\ell_1(z_i;\mu_2,\sigma_2, \xi_2)}\right]. 
\end{eqnarray}
Consequently, the maximum-likelihood estimation of $\boldsymbol\theta$ can be obtained by maximizing the log-likelihood function in Eq.\eqref{log-likelihood-acc-pmax}. 
Similarly, if $Z \sim \underline{H}$, an accelerated $p$-min stable distribution, one may find its maximum likelihood estimators. Note that the maximum likelihood estimations involved in the GEV models don't exist when $\xi \le -1$. In the following, we will establish the consistency of mle in Theorem \ref{Theoremconsistency} below when both extreme value index is greater than $-1$.

\textbf{Consistency and asymptotic normality of the maximum likelihood estimation}. 
Note that the mle of the three parameters of $(\xi,\mu,\sigma)$ for the $p$-max stable distribution and $p$-min stable distribution  $H_1(x;\mu, \sigma, \xi)=G_\xi(\log x;\mu,\sigma)$ and $\underline H_1(x) = 1- H_1(-x; \mu, \sigma, \xi)$ exists with consistency for $\xi>-1$, and asymptotic normality for $\xi>-1/2$ \citep{Dombry2015,SmithRichardL.1985Mlei}. In the following, we will consider the properties of the mle of the parameters involved in the accelerated $l$-min stable distribution, when it exists. This will indicate the same results for the classes of accelerated $l$-max, $p$-max and $p$-min stable distributions by corresponding transformations, e.g., as shown in Corollary \ref{Cor-1}. {The main idea is to apply \citet[Theorems 1 and 3]{SmithRichardL.1985Mlei} to show its consistency  and asymptotic normality in Theorems  \ref{Theoremconsistency} and  \ref{Theoremasymptotic normality}, respectively.} Recall that the risk $X$ considered in \citet{SmithRichardL.1985Mlei} has  probability density function as 
\begin{equation}\label{Eq: Smith}
    f(x;\theta, \boldsymbol{\phi})=(x-\theta)^{\alpha -1}g(x-\theta;\boldsymbol{\phi}), \quad\quad \theta<x<\infty,
\end{equation}

where $\alpha>1, \theta$ and $\bm{\phi}=(\phi_1,\phi_2,\cdots,\phi_p) \in \Phi \subset \mathbb{R}^{p}$ are unknown parameters, which might a function of $\alpha$, and $g(x; \bm{\phi})\to \alpha c(\bm{\phi})>0$ as $x \to 0$. 

Note that the accelerated $l$-min stable distribution is given by

$$\underline H(x; \mu_1, \sigma_1, \xi_1, \mu_2, \sigma_2, \xi_2) = 1-G_{\xi_1}(-x; \mu_1, \sigma_1)G_{\xi_2}(-x; \mu_2, \sigma_2).$$ 

In what follows, we first re-parameterize the accelerated $l$-min stable density function $\underline h=\underline H'$ as in Eq.\eqref{Eq: Smith}. Denote by $\alpha_1=-1/\xi_1, \alpha_2=-1/\xi_2$ and assume $\alpha_1,\alpha_2>0$ 
and $\theta_1=-\mu_1+\sigma_1/\xi_1\le \theta_2=-\mu_2+\sigma_2/\xi_2$ without loss of generality. Hence, its density is given by
\begin{eqnarray*}
&&\underline h(x;\theta_1,\sigma_1,\alpha_1,\theta_2,
    \sigma_2,\alpha_2)\\
    &&\quad =\begin{cases} 
    \frac{1}{\sigma _1} \left ( \frac{x-\theta_1}{\sigma_1\alpha_1}  \right )^{\alpha_1-1} \exp\left[-\left ( \frac{x-\theta _1}{\sigma _1\alpha _1}  \right ) ^{\alpha_1}\right ], &\theta_1 < x\le \theta_2,\\
\left [ \frac{1}{\sigma _1} \left ( \frac{x-\theta_1}{\sigma_1\alpha_1}  \right )^{\alpha_1-1}+ \frac{1}{\sigma _2} \left ( \frac{x-\theta_2}{\sigma_2\alpha_2}  \right )^{\alpha_2-1}  \right ]\exp\left[-\left ( \frac{x-\theta _1}{\sigma _1\alpha _1}  \right ) ^{\alpha_1}-\left ( \frac{x-\theta _2}{\sigma _2\alpha _2}  \right ) ^{\alpha_2}\right],&x>\theta _2. 
\end{cases}
\end{eqnarray*}
Below, we consider the case with $\theta_1=\theta_2=\theta$ and $1<\alpha_1 \le
\alpha_2$\ (recall $\alpha_1$ and $\alpha_2$ are exchangeable), and thus the pdf $f$ can be written of form Eq.\eqref{Eq: Smith}, i.e., 
\begin{eqnarray}
    \label{acc-l-min-pdf}
\underline h(x;\theta,\sigma_1,\alpha_1,\sigma_2,\alpha_2)&=&  
 (x-\theta)^{\alpha_1-1}\left [ \frac{1}{\sigma _1} \left ( \frac{1}{\sigma_1\alpha_1}  \right )^{\alpha_1-1}+ \frac{1}{\sigma _2} \left (\frac{1}{\sigma_2\alpha_2}  \right )^{\alpha_2-1} (x-\theta)^{\alpha_2-\alpha_1} \right ] \notag \\
 && \times \exp\left[-\left ( \frac{x-\theta }{\sigma _1\alpha _1}  \right ) ^{\alpha_1}-\left ( \frac{x-\theta }{\sigma _2\alpha _2}  \right ) ^{\alpha_2}\right]\notag\\
 &=:& (x-\theta)^{\alpha -1} g(x-\theta; \bm{\phi}),\qquad x>\theta
\end{eqnarray}
 with $\alpha=\alpha_1$ and $\bm{\phi}=(\sigma_1,\alpha_1,\sigma_2,\alpha_2) \in \Phi := R^+\times (1,\infty) \times R^+\times (1,\infty)$ and 
\begin{eqnarray} \label{g-acc-l-min}
    g(x;\bm{\phi})&=&\left [ \frac{1}{\sigma _1} \left ( \frac{1}{\sigma_1\alpha_1}  \right )^{\alpha_1-1}+ \frac{1}{\sigma _2} \left ( \frac{1}{\sigma_2\alpha_2}  \right )^{\alpha_2-1}  x^{\alpha_2-\alpha_1} \right ] \notag\\
    &&\times \exp\left[-\left ( \frac{x}{\sigma _1\alpha _1}  \right ) ^{\alpha_1}-\left ( \frac{x}{\sigma _2\alpha _2}  \right ) ^{\alpha_2}\right], \quad x>0. 
\end{eqnarray}

\begin{theorem}\label{Theoremconsistency}
   Let $X_1,\ldots, X_n$ be a random sample from an accelerated $l$-min stable distribution with density $\underline h(x;\theta,\bm\phi)$ given by  Eq.\eqref{acc-l-min-pdf} 
   and $\xi_1,\xi_2<0$. If  $\alpha = \alpha_1 = -1/\xi_1, \alpha_2 = -1/\xi_2$  satisfy that (i) $\alpha>1$, and $\alpha_1 = \alpha_2$ or $\alpha_1 < \alpha_2 -1$; or (ii) $\alpha>2$ and $\alpha_1 < \alpha_2 -2$,
 then there exist  a sequence $(\widehat{\theta}_n, \widehat{\bm{\phi}}_n)$ as the solution of 
the likelihood equations (recall $\alpha_1 = \alpha, \bm \phi = (\sigma_1,\alpha_1,\sigma_2,\alpha_2)$): 
\begin{equation*}
    \frac{\partial \ell_n \left (\theta,\bm{\phi}  \right ) }{\partial \theta}=0, \quad \frac{\partial \ell_n \left (\theta,\bm{\phi} \right ) }{\partial \phi_i}=0,\quad i=1,\ldots,4
\end{equation*}
  with the log-likelihood function 
   $\ell_n(\theta,\bm\phi)$  given by
\begin{equation*}
\ell_n(\theta,\bm{\phi})=\sum_{i=1}^{n} \log \underline h\left ( X_i;\theta ,\bm{\phi} \right ).
\end{equation*}
We have $\widehat{\theta}_n$ and $\widehat{\bm{\phi}}_n$ are the consistent estimators of $\theta$ and $\bm\phi$, respectively. 
\end{theorem}

Given $\underline h(x;\alpha, \bm\phi)$ in Eq.\eqref{acc-l-min-pdf}, define  by $\bm M$ the matrix as $(m_{ij}(\bm{\phi}))_{i,j=0,\cdots,4}$ for $\alpha>2$, and $(m_{ij}(\bm{\phi}))_{i,j=1,\cdots,4}$ for $1<\alpha \leq 2$. Here
\begin{equation*}
    m_{ij}(\bm{\phi})=-\EE_{\bm{\phi}} \left \{\frac{\partial^2}{\partial \phi_i\partial \phi_j}\log \underline h(X;0,\bm{\phi} )  \right \}, i,j=1,\ldots,4,
\end{equation*}
\begin{equation*}
    m_{0i}(\bm{\phi})=m_{i0}(\bm{\phi})=\EE_{\bm{\phi}} \left \{\frac{\partial^2}{\partial x\partial \phi_i}\log \underline h(X;0,\bm{\phi} )  \right \}, i=1,\ldots,4
\end{equation*}
and
\begin{equation*}
    m_{00}(\bm{\phi})=-\EE_{\bm{\phi}} \left \{\frac{\partial^2}{\partial x^2}\log \underline h(X;0,\bm{\phi} )  \right \}.
\end{equation*}

\begin{theorem}\label{Theoremasymptotic normality}
    Under the assumptions of Theorem \ref{Theoremconsistency}, let $(\widehat{\theta}_n,\widehat{\bm{\phi}}_n)$ denote the maximum likelihood estimators of the parameters involved in the accelerated $l$-min stable distribution with pdf defined by Eq.\eqref{acc-l-min-pdf}. 
    \begin{itemize}
        \item[(i).] If $\alpha = \alpha_1 >2$ and $\alpha_1<\alpha_2-2$, then $\sqrt n(\widehat{\theta}_n-\theta, \widehat{\bm\phi}_n-\bm\phi)$ converges in distribution to a Gaussian random vector with mean ${\bf 0}$ and variance-covariance matrix $\bm{M}^{-1}$.
        \item[(ii).] If $\alpha= \alpha_1 =2$ and $\alpha_1<\alpha_2-1$, then  $\left \{ \sqrt{(\sigma_1\alpha_1)^{-\alpha_1}n\log n}(\widehat{\theta}_n-\theta),\sqrt{n}(\widehat{\bm\phi}_n-\bm\phi) \right \} $ converges in distribution to a Gaussian random vector with mean zero and variance-covariance matrix of the form 
        \begin{equation*}
            \begin{bmatrix}
             1&0 \\
             0&\bm{M}^{-1} 
            \end{bmatrix}.
        \end{equation*}
       \end{itemize}
\end{theorem}
We will only provide the proof of Theorem \ref{Theoremconsistency} in Appendix \ref{Appendix C}, since Theorem \ref{Theoremasymptotic normality} follows immediately by Theorem \ref{Theoremconsistency} and \citet[Theorem 3 (i) and (ii)]{SmithRichardL.1985Mlei}.

\begin{remark}
 (i) The pseudo maximum likelihood estimators of $(\theta, \bm\phi)$ obtained in Theorem \ref{Theoremconsistency} is essentially the global mle due to its positive definite information matrix for $-1/2 < \xi_1, \xi_2<0$.  \\
 (ii) In addition, similar to the GEV,  the maximum likelihood estimators exist also with consistency and asymptotic normality for the parameters involved in the accelerated $l$-min stable distribution with  both $\xi_1$ and $\xi_2$ being positive and $\theta_1 := -\mu_1+\sigma_1/\xi_1 = -\mu_2+\sigma_2/\xi_2 =: \theta_2$.
\end{remark}

\section{Examples}
\label{sec: Examples}

This section provides illustrating examples to demonstrate the main findings of the paper.  We will see how the sample generating process and interplay between competing risks result in the three limit types, namely the accelerated $p$-max stable distributions, $p$-max stable distributions and left-truncated $p$-max stable distribution according to Theorem \ref{TheoremMaxMax}. 
\begin{example}[Log-Fr{\'e}chet and log-Fr{\'e}chet]\label{FandF}
 Let  $F_1(x)=H _{1,\alpha_1} (x)$
 and $F_2(x)=H_{1,\alpha_2}(x)
 $ with $\alpha_1>\alpha_2>0$, which are called log-Fr\'echet distributions since $\widetilde F_j(x) := F_j(\e^x)$ is Fr\'echet distribution. Hence, both $F_j$'s are super heavy tailed distributions, it does not belong to any $l$-MDAs. Meanwhile, $F_j$ is $p$-max stable distribution, and $F_j\in  D_{p}\left ( H_{1,\alpha _{j} } \right )$. Indeed, Eq.\eqref{ComponentPowerLimit} holds with 
  (cf. \citet[Theorem 2.1]{Mohan1993})
\begin{eqnarray}\label{p-Norm-logFrechet}
    \alpha_{1,n_1}= \alpha_{2,n_2} =1, \quad \beta_{j,n_j} =
n_j^{-1/\alpha_j}, \quad H^{(j)}(x) = H_{1,\alpha_j}(x)\quad j=1,2.
\end{eqnarray}
Therefore, we  have (recall Eq.\eqref{NormingMaxMax})  

$$\alpha_n=\alpha_{1,n_1}\left(\frac{1}{\alpha_{2,n_2}}\right)^\frac{\beta_{1,n_1}}{\beta_{2,n_2}} = 1,\quad \beta_n=\frac{\beta_{1,n_1}}{\beta_{2,n_2}} = 
\frac{n_2^{1/\alpha_2}}{n_1^{1/\alpha_1}}$$

implying thus $A=1.$ Note that the limit $\beta_n$ depends on both sample sizes and tail parameters. Concerning the sample size, we consider the following two cases: \underline{(I) $n_2=cn_1$ with some $c>0$} and \underline{(II) $n_2=an_1^c$ with some $a,c>0$}.

\underline{Case I: $n_2=cn_1$ with some $c>0$}. Note that $\alpha_1>\alpha_2>0$, and thus 
$\beta_n \sim c^{1/\alpha_2} n_1^{1/\alpha_2-1/\alpha_1} \to B$ with $B =\infty$. Noting further both supports of $H^{(j)}$ are $(1, \infty)$. Therefore, $x_n = \alpha_n x^{\beta_n} \to \infty$ for all $x>1$,  implying that $H^{(1)}(x_n)=H_{1,\alpha_1}(x_n) \to 1$, and thus Theorem \ref{TheoremMaxMax}(ii) holds, i.e., 
\begin{equation*}
   \mathbb{P} (\alpha_{2,n_2}| M_n |^{\beta_{2,n_2}}\sign(M_n) \le x) \toweak H_{1,\alpha_2}(x).
\end{equation*}
This indicates that the maxima from the heavier distribution $F_2$ overwhelms the lighter one when their sample size is of proportion.

\underline{Case II: $n_2=an_1^c$ with some $a,c>0$}. We have $\beta_n \sim n_1^{{c}/{\alpha_2}-1/{\alpha_1}}$. The limit of $\beta_n$ is $B=\infty, 1,0$ provided that $c>, =, <{\alpha_2}/{\alpha_1}$ accordingly. Hence, Theorem \ref{TheoremMaxMax}(i)  holds if $ n_2 = b^{\alpha_2}n_1^{\alpha_2/\alpha_1}$, $b>0$, i.e., 
\begin{equation*}
    \mathbb{P}(\alpha_{2,n_2}|M_n|^{\beta_{2,n_2}}\sign(M_n)\le x ) \toweak H_{1,\alpha_1}(|x|^b \sign(x)) H_{1,\alpha_2}(x).
\end{equation*}
We see that the competing maxima converge to the accelerated $p$-max stable distribution, which is not of any type of power maximum domain attraction. Additionally, we see that, although $F_2$ is heavier than $F_1$, the sample maxima from $F_1$ over a larger length is comparable, even overwhelms $F_2$ if $n_1$ tends to infinity fairly faster as shown below.
 While Theorem \ref{TheoremMaxMax}(ii) holds for $c> {\alpha_2}/{\alpha_1}$, following similar arguments of Case I. 
Noting that for $0<c<\alpha_2/\alpha_1$, $1/\beta_n \to \infty$, we have 
\begin{equation*}
    \mathbb{P} (\alpha_{1,n_1}| M_n |^{\beta_{1,n_1}} \sign(M_n)\le x) \toweak H_{1,\alpha_1}(x).
\end{equation*}
\end{example}

\begin{example}[Log polynomial and log polynomial]\label{WandW}
Log polynomial distributions $F_j(x) = F(\log x; \alpha_j)$ where $F(x; \alpha)=1-(\gamma(F)-x)^\alpha, \gamma(F)-1<x<\gamma(F)$ is a polynomial distribution with parameter $\alpha>0$. 
The choices of power normalization constants can be referred to Table \ref{Table: CommonDistribution}.


In particular, we consider the following two cases in power normalization: \underline{(I) $n_2=cn_1$ with some $c>0$} and\\
\underline{(II) $n_2=an_1^c$ with some $a, c>0$}.

\underline{Case I: $n_2=cn_1$ with some $c>0$}. 
We have  
\begin{equation*}
    \mathbb{P}(\alpha_{2,n_2}|M_n|^{\beta_{2,n_2}}\sign(M_n) \le x ) \toweak H_{2,\alpha_2}(x).
\end{equation*}
\underline{Case II: $n_2=an_1^c$ with some $a,c>0$}. 
For $ n_2 = b^{-\alpha_2}n_1^{\alpha_2/\alpha_1}$, $b>0$, $\beta_n \to B=b$, Theorem \ref{TheoremMaxMax} (i) holds with
 
\begin{equation*}
    \mathbb{P}(\alpha_{2,n_2}|M_n|^{\beta_{2,n_2}}\sign(M_n) \le x ) \toweak H_{2,\alpha_1}(|x|^b\sign(x)) H_{2,\alpha_2}(x).
\end{equation*}
Note that for $0<c<{\alpha_2}/{\alpha_1}, \beta_n  \to \infty$, we have {following similar arguments for $1/B=0$ as for Case I}
{\begin{equation*}
    \mathbb{P}(\alpha_{1,n_1}|M_n|^{\beta_{1,n_1}}\sign(M_n) \le x ) \toweak H_{2,\alpha_1}(x).
\end{equation*}
}
In linear normalization, we can choose linear normalization constants from Table \ref{Table: CommonDistribution} to have
\begin{equation*}
    \pk{a_{j,n_j}(M_{j,n_j}-b_{j,n_j})\leq x}\toweak \Psi_{\alpha_j}(x),\quad j=1,2.
\end{equation*}
Recalling \citet[Theorem 2.1]{cao2021}, we will have the following claims for the max of max under linear normalization following 
similar arguments as above: 
\begin{itemize}
    \item For $n_2 = cn_1, c>0$ or $n_2 = n_1^c$ with $c>\alpha_2/\alpha_1$, it follows by 
    $a_n\to 0$ that 
    \begin{equation*}
    \pk{a_{2,n_2}(M_n-b_{2,n_2})\leq x}\toweak \Psi_{\alpha_2}(x).
\end{equation*}
\item For $n_2 = n_1^c$ with $c = \alpha_2/\alpha_1$, it follows by $a_n\to1$ that
    \begin{equation*}
    \pk{a_{2,n_2}(M_n-b_{2,n_2})\leq x}\toweak \Psi_{\alpha_1}(x)\Psi_{\alpha_2}(x).
\end{equation*}
\item For $n_2 = n_1^c$ with $c<\alpha_2/\alpha_1$, it follows by $1/a_n\to 0$ that
\begin{equation*}
     \pk{a_{1,n_1}(M_n-b_{1,n_1})\leq x}\toweak \Psi_{\alpha_1}(x).
\end{equation*}

\end{itemize}

\COM{In particular, we consider the following two cases \underline{(I) $n_2=c n_1$ with some $c>0$} and \underline{(II) $n_2 = n_1^c$ with some $c>0$}.

\underline{Case I: $n_2=c n_1$ with some $c>0$}. Note that $\alpha_1>\alpha_2$, and thus $a_n \sim c^{-1/\alpha_2}(n_1)^{1/\alpha_1-1/\alpha_2} \to 0$, thus we have
\begin{equation*}
    \pk{a_{2,n_2}(M_n-b_{2,n_2})\leq x}\toweak \Psi_{\alpha_2}(x).
\end{equation*}
\underline{Case II: $n_2 = n_1^c$}. We have $a_n \sim n_2^{1/\alpha_1-c/\alpha_2}$. Hence for $c>\alpha_2/\alpha_1$, $a_n \to 0$, we have
\begin{equation*}
    \pk{a_{2,n_2}(M_n-b_{2,n_2})\leq x}\toweak \Psi_{\alpha_2}(x).
\end{equation*}
For $c=\alpha_2/\alpha_1$, $a_n \to 1$, we have
\begin{equation*}
    \pk{a_{2,n_2}(M_n-b_{2,n_2})\leq x}\toweak \Psi_{\alpha_1}(x)\Psi_{\alpha_2}(x).
\end{equation*}
Note that for $0<c<\alpha_2/\alpha_1$, $a_n \to \infty$, we have 
\begin{equation*}
     \pk{a_{1,n_1}(M_n-b_{1,n_1})\leq x}\toweak \Psi_{\alpha_1}(x).
\end{equation*}
}
\end{example}

\begin{example}[Pareto and log-Fr{\'e}chet]
\label{ex33}
Let 
$F_{1}(x) = 1- x^{-\alpha_1}, x>1$ be a Pareto distribution with parameter $\alpha_1>0$ 
and $F _{2}\left ( x \right )=  H_{1,\alpha_2}(x)
$. We have $F_1\in D_l(\Phi_1)$ and $F_1\in D_p(H_5)$. The choices of normalization constants can be referred to Table \ref{Table: CommonDistribution}.

Therefore, we have $\alpha_n = \alpha_{1,n_1} = 1/n_1 \to A = 0$ and $\beta_n = \alpha_1 n_2^{1/\alpha_2} \to B=\infty$. Given the support of $H^{(2)}$ is  $(1,\infty)$, we need to check if $x_n = \alpha_n x^{\beta_n}$ tends to infinity for any $x>1$, i.e.,
$$\log x_n =   \alpha_1 n_2^{1/\alpha_2}\log x - \log n_1 \to\infty.$$
Clearly, for $n_2=cn_1, n_1^c$ with some $c>0$ or $n_2 = (\log n_1)^c$ for some $c>\alpha_2$, we have $x_n\to\infty$. Thus Theorem \ref{TheoremMaxMax}(ii) holds. For $n_2 = (\log n_1/\alpha_1)^{\alpha_2}, x_n \to \I{x > \e}\infty$ with $\e \in (1,\infty)$, the support of $H_{1,\alpha_2}$, thus we have (cf. Theorem \ref{TheoremMaxMax}(iii))
\begin{equation*}
    \mathbb{P}(\alpha_{2,n_2}|M_n|^{\beta_{2,n_2}} \sign(M_n)\le x ) \to H_{1,\alpha_2}(x)\I{x > \e}.
\end{equation*}
For $n_2 = (\log n_1)^c, 0<c<\alpha_2$, we have 
\begin{equation*}
    \mathbb{P}(\alpha_{1,n_1}|M_n|^{\beta_{1,n_1}} \sign(M_n)\le x ) \to H_{5}(x).
\end{equation*}
\end{example}

\COM{
\begin{example}[Normal and Fr{\'e}chet]\label{NandF}
Let $F_{1}(x)=\Phi\left (x\right )$ be the standard normal distribution and $F_{2}\left(x\right) = \Phi_{\alpha}(x)$ be the Fr\'echet distribution. We have $F_1\in D_l(\Lambda)$ and $F_2\in D_l(\Phi_{\alpha})$. In addition, both belong to $p$-MDA such that Eq.\eqref{ComponentPowerLimit} holds with (cf. \citet[Theorems 2.1]{Mohan1993}, \cite{ref8} and \cite{ref9})
\begin{eqnarray*}
    \alpha_{1,n_1} &=& \left(\left(2 \log n_1\right)^{1 / 2}-\frac{1}{2}\left(2 \log n_1\right)^{-1 / 2}\left(\log \log n_1+\log 4 \pi\right)\right)^{-\beta_{1,n_1}} \approx (2\log n_1)^{-\log n_1},\\
    \beta_{1,n_1} &=& 2 \log n_1 -\frac{1}{2}\left(\log \log n_1+\log 4 \pi \right) \approx 2\log n_1, \\
   \alpha_{2,n_2} &=& 1/{n_2}, \quad \beta_{2,n_2} = \alpha,\quad H^{(1)} = H^{(2)}= H_{5}. 
\end{eqnarray*}
 \COM{{Therefore, we have $\alpha_n \sim \left(2\log n_1 \right)^{-\log n_1}$ and $\beta_n \sim (2\log n_1) n_2^{1/\alpha} $. 
Similar arguments for Example \ref{ex3}, we analyse the cases that 
$$\log x_n =   (\log n_1)\big[2n_2^{1/\alpha}\log x - \log\log n_1 - \log 2 \big]\to\infty\quad \mbox{for}\ x>1.$$
Clearly, for $n_2=cn_1, n_1^c, (\log n_1)^c$ with some $c>0$ or $n_2 = (\log \log n_1)^c$ for some $c>\alpha$, we have $x_n\to\infty$. Thus Theorem \ref{TheoremMaxMax}(ii) holds. For $n_2 = (\log \log n_1)^\alpha, x_n \to \I{x>{\sqrt e}}\infty$, thus we have 
\begin{equation*}
    \mathbb{P}(\alpha_{2,n_2}|M_n|^{\beta_{2,n_2}} \sign(M_n)\le x ) \toweak H_{1,\alpha}(x)\I{x>{\sqrt e}}.
\end{equation*}
For $n_2 = (\log \log n_1)^c, 0<c<\alpha$, we have 
\begin{equation*}
    \mathbb{P}(\alpha_{1,n_1}|M_n|^{\beta_{1,n_1}} \sign(M_n)\le x ) \toweak H_{5}(x).
\end{equation*}
} }
Therefore, we have $\alpha_n \sim \left[n_2^{2/\alpha}/(2\log n_1) \right]^{\log n_1}$ and $\beta_n \sim (2\log n_1)/\alpha $. Given the support for both $H^{(1)}$ and $H^{(2)}$ is $(0,\infty)$, we analyse the cases that 
$$x_n =  \left[ (n_2x)^{2/\alpha}/2\log n_1 \right]^{\log n_1} \to\infty\quad \mbox{for}\ x>0.$$
Clearly, for $n_2=cn_1, n_1^c$ with some $c>0$ or $n_2 = (\log n_1)^c$ for some $c>\alpha/2$, we have $x_n\to\infty$. Thus Theorem \ref{TheoremMaxMax}(ii) holds. For $n_2 = (\log n_1)^{\alpha/2}, x_n \to \I{x> 2^{\alpha/2}}\infty$ with $2^{\alpha/2}\in(1,\infty)$, the support of $H_{1, \alpha}$, thus we have (see Theorem \ref{TheoremMaxMax}(iii)) 
\begin{equation*}
    \mathbb{P}(\alpha_{2,n_2}|M_n|^{\beta_{2,n_2}}\sign(M_n)\le x ) \to H_{5}(x)\I{x>2^{\alpha/2}}.
\end{equation*}
For $n_2 = (\log n_1)^c, 0<c<\alpha/2$, we have 
\begin{equation*}
    \mathbb{P}(\alpha_{1,n_1}|M_n|^{\beta_{1,n_1}}\sign(M_n)\le x ) \to H_{5}(x).
\end{equation*}

To compare with the linear normalization maxima of maxima, the limit distribution exists only with (recall $a_{1,n_1} =\sqrt{2\log n_1}, b_{1,n_1}=\sqrt{2\log n_1} -(1/2)(2\log n_1)^{-1/2}\left(\log \log n_1+\log 4 \pi \right)$)
\begin{eqnarray*}
    \pk{a_{1,n_1}(M_n - b_{1,n_1}) \le x} \to \Lambda(x)
\end{eqnarray*}
as $n_2^{2/\alpha}/\log n_1$ tends to zero, following similar arguments as for \cite[Example 2.7]{cao2021}. We see that the limit distributions of power normalized max of max are fairly wider than the linear normalization.

\end{example}
}

\COM{
\begin{example}[Skew-Normal and log-Fr{\'e}chet]
Let $F_{1}\left ( x\right )=\phi \left ( x \right ) \Phi \left ( \lambda x \right )$ be the skew-normal distribution introduced in \cite{ref10}, and $F_2(x)=H_{1,\alpha}(x)$. Note that $\phi \left ( x \right )$ is the pdf of the standard normal distribution and $\Phi\left ( x \right )$ is the cumulative distribution function (cdf) of standard normal distribution, 
{and $\lambda$ is} the shape parameter $\lambda \in \mathbf{R}$. According to \cite{ref8}, \cite{ref9} and \citet[Theorems 2.1]{Mohan1993}, $F_1\in D_l(\Lambda)$ with $x_{F_1} = \infty$ and $F_1 \in D_p(H_5)$ for $\lambda>0$ and $F_2\in D_p(H_{1,\alpha})$, with the normalization constants 
\begin{eqnarray*}
{a_{1, n_1}} &=& b_{1, n_1} = \gamma_{1, n_1}(F_\lambda) = \sqrt{2\log n_1}, \ G(x) = \expon(-\expon(-x)),\\
    \alpha_{1,n_1} &= & \kai{[\gamma_{1, n_1}(F_\lambda)]
   {^{-\gamma^2_{n_1}(F_\lambda)}} 
    }  \approx (2\log n_1)^{\kai{-\log n_1}}, \ \beta_{1,n_1} = \kai{\gamma_{1, n_1}^2(F_\lambda)} \approx 2\log n_1, \\
    \alpha_{2,n_2} &=& 1, \quad \beta_{2,n_2} = n_2^{-\frac{1}{\alpha}},\quad H^{(2)} = H_{1,\alpha}, \ H^{(1)}= H_{5}. 
\end{eqnarray*}
Clearly, for $n_2=cn_1, n_1^c, (\log n_1)^c$ with some $c>0$ or $n_2 = (\log \log n_1)^c$ for some $c>\alpha$, we have $x_n\to\infty$. Thus Theorem \ref{TheoremMaxMax}(ii) holds. For $n_2 = (\log \log n_1)^\alpha, x_n \to \I{x> \kai{\sqrt{e}}}\infty$, thus we have 
\begin{equation*}
    \mathbb{P}(\alpha_{2,n_2}|M_n|^{\beta_{2,n_2}}\sign(M_n)\le x ) \toweak H_{1,\alpha}(x)\I{x> \kai{\sqrt{e}}}.
\end{equation*}
For $n_2 = (\log \log n_1)^c, 0<c<\alpha$, we have 
\begin{equation*}
    \mathbb{P}(\alpha_{1,n_1}|M_n|^{\beta_{1,n_1}}\sign(M_n)\le x ) \toweak H_{5}(x).
\end{equation*}
\end{example}
}

\begin{example}[Uniform and Uniform]\label{UandU}
Let $F_1$ and $F_2$ be $U(l_1,u_1)$ and $U(l_2,u_2)$ with both $u_2>u_1>0$, two uniform distributions, respectively. Thus, both $F_1, F_2$ belong to $D_l(\Psi_\alpha)$ with Eqs.\eqref{l-MDA} and \eqref{p-MDA} satisfied with linear and power normalization constants in Table \ref{Table: CommonDistribution}.

For linear normalization, we have
\begin{equation*}
   \mathbb{P} (a_{2,n_2} (M_n - b_{2,n_2}) \le x) \toweak \Psi_1(x).
\end{equation*}
And for power normalization, we have  
\begin{equation*}
   \mathbb{P} (\alpha_{2,n_2}| M_n |^{\beta_{2,n_2}}\sign(M_n) \le x) \toweak H_{2,1}(x).
\end{equation*}
\end{example}

\begin{example}[General error and log-Fr{\'e}chet]\label{GandF}

Let $F_{1}$ be the general error distribution introduced by \citet{Chen2012}, and $F_{2}\left(x\right) = H_{1,\alpha }\left( x\right)$. Note that the pdf of $F_{1}$ is 
given  by
\begin{equation*}
    f_{1}(x; \nu) =\frac{2^{-2/\nu}\expon { -\frac{1}{2}  \left | \frac{x}{\lambda }  \right |^{\nu}}   }{\lambda 2^{1+1/\nu }\Gamma \left ({\nu^{-2/\nu}}  \right ) }, \quad x\in\mathbb R,  
\end{equation*}
where $\nu>0$ and $\lambda =\sqrt{{2^{-2/\nu}}\Gamma \left ( 1/{\nu}  \right )\Gamma \left ( 3/{\nu}  \right )}, \Gamma(\cdot)$ is the gamma function. 
Both $F_1$ and $F_2$ belong to $p$-MDA such that Eq.\eqref{ComponentPowerLimit} holds with (cf. \citet[Theorem 2.1]{Mohan1993} and \citet{Chen2012}) normalization constants in Table \ref{Table: CommonDistribution}.

Therefore, for $n_2=cn_1, n_1^c, {(\log n_1)^c}$ with some $c>0$ or $n_2 = (\log \log n_1)^c$ for some $c>\alpha$, we have $x_n\to\infty$. Thus Theorem \ref{TheoremMaxMax}(ii) holds. For $n_2 = (\log \log n_1)^\alpha, x_n \to \I{x>{\e^{1/\nu}}}\infty$ with $\e^{1/\nu}\in(1,\infty)$, the support of $H_{1, \alpha}$, thus we have (cf. Theorem \ref{TheoremMaxMax}(iii))
\begin{equation*}
    \mathbb{P}(\alpha_{2,n_2}|M_n|^{\beta_{2,n_2}} \sign(M_n)\le x ) \to H_{1,\alpha}(x)\I{x> \e^{1/\nu}}.
\end{equation*}
For $n_2 = (\log \log n_1)^c, 0<c<\alpha$, we have 
\begin{equation*}
    \mathbb{P}(\alpha_{1,n_1}|M_n|^{\beta_{1,n_1}} \sign(M_n)\le x ) \to H_{5}(x).
\end{equation*}

\end{example}

\begin{example}[Pareto and Pareto]\label{ex3.8}
Let  $F_j(x) = 1- x^{-\alpha_j}, x>1$ be Pareto distributions or Fr\'echet distribution with parameters $\alpha_j>0, \, j=1,2$. We see that $F_j \in D_l(\Phi_{\alpha_j})$ and $F_j \in D_p(H_5)$. Indeed,  Eqs.\eqref{l-MDA} and \eqref{ComponentPowerLimit} hold with (cf. \citet[Theorems 2.5, Theorem 3.1(a)]{Mohan1993}) power normalization constants in Table \ref{Table: CommonDistribution}.

Below, we discuss the limit distribution according to the limit of $\alpha_n$.\\
\underline{Case I}. For $n_2 = cn_1^{\alpha_2/\alpha_1}$, Theorem \ref{TheoremMaxMax} (i) holds with $A=c, B=\alpha_1/\alpha_2$, i.e., 
\begin{equation*}
    \mathbb{P}(\alpha_{2,n_2}|M_n|^{\beta_{2,n_2}} \sign(M_n)\le x ) \toweak H_{5}(c|x|^{\alpha_1/\alpha_2}\sign(x)) H_{5}(x). 
\end{equation*}
\underline{Case II}. For $n_2/n_1^{\alpha_2/\alpha_1} \to\infty$, e.g., $n_2 = cn_1$, Theorem \ref{TheoremMaxMax} (ii(b)) holds with $A=\infty, B=\alpha_1/\alpha_2$, i.e., 
\begin{equation*}
    \mathbb{P}(\alpha_{2,n_2}|M_n|^{\beta_{2,n_2}}  \sign(M_n)\le x ) \toweak H_{5}(x). 
\end{equation*}
\underline{Case III}. For $n_2=o(n_1^{\alpha_2/\alpha_1})$, e.g., $n_2 = n_1^c$ with $0<c<\alpha_2/\alpha_1$, we have 
\begin{equation*}
    \mathbb{P}(\alpha_{1,n_1}|M_n|^{\beta_{1,n_1}} \sign(M_n)\le x ) \toweak H_{5}(x). 
\end{equation*}

In linear normalization, we can choose linear normalization constants to have
\begin{equation*}
    \pk{a_{j,n_j}(M_{j,n_j}-b_{j,n_j})\leq x}\toweak \Phi_{\alpha_j}(x),\quad j=1,2.
\end{equation*}
Recalling \citet[Theorem 2.1]{cao2021}, we will have the following claims for the max of max under linear normalization with
similar arguments as above: 
\begin{itemize}
    \item For $n_2 = cn_1, c>0$ or $n_2 = n_1^c$ with $c>\alpha_2/\alpha_1$, it follows \cL{by} 
    $a_n\to0$ that 
    \begin{equation*}
    \pk{a_{2,n_2}(M_n-b_{2,n_2})\leq x}\toweak \Phi_{\alpha_2}(x).
\end{equation*}
\item For $n_2 = n_1^c$ with $c = \alpha_2/\alpha_1$, it follows by $a_n\to1$ that
    \begin{equation*}
    \pk{a_{2,n_2}(M_n-b_{2,n_2})\leq x}\toweak \Phi_{\alpha_1}(x)\Phi_{\alpha_2}(x).
\end{equation*}
\item For $n_2 = n_1^c$ with $c<\alpha_2/\alpha_1$, it follows by $1/a_n\to 0$ that
\begin{equation*}
     \pk{a_{1,n_1}(M_n-b_{1,n_1})\leq x}\toweak \Phi_{\alpha_1}(x).
\end{equation*}
\end{itemize}
\end{example}

\begin{example}[Polynomial growth 
and Fr{\'e}chet]\label{ex3.9}
\label{PandF}
Let  $F_1(x)=1-(\gamma(F_1)-x)^{\alpha_1}, \gamma(F_1)-1<x<\gamma(F_1)$ with $\gamma(F_1)>0$
and $F_2(x)=\Phi_{\alpha_2}(x)$. We have $F_1\in D_{l}\left (\Psi_{\alpha_1} \right)$ and $F_2 \in D_{l}\left (\Phi_{\alpha_2} \right )$. 
In addition, $F_1\in D_p(H_{2, \alpha_1})$ and $F_2\in D_p(H_5)$ 
 such that Eq.\eqref{ComponentPowerLimit} holds with (cf. \citet[Theorems 2.2 and 2.5]{Mohan1993}) normalization constants provided in Table \ref{Table: CommonDistribution}.
We have 
$$\log x_n = \log \alpha_n +{\beta_n} \log x \to \infty$$ holds for any $x>0$. 
Consequently,  the claim of  Theorem \ref{TheoremMaxMax}(ii) follows.

Note that if we set the right-endpoint $\gamma(F_1)<0$ for the polynomial growth distribution, then $F_1\in D_p(H_{4, \alpha_1})$ and Theorem \ref{TheoremMaxMax}(ii(a)) follows straightforwardly.
\end{example}
To summarize, the power-normalized extremes behave quite differently. In general, the limit types are determined basically by the tail heaviness of the potential risks and the convergence of the sample maxima to their right endpoints. For the convenience of further simulations, Table \ref{Table: CommonDistribution} summarizes common examples in one/both max-domain attractions along with their normalization constants.

\begin{table}[H]
    \centering
    \caption{{Common distributions in $l$-max domain attractions and $p$-max domain attractions with its normalization constants satisfying Eqs.\eqref{l-MDA} and \eqref{p-MDA}.}}
    \resizebox{\linewidth}{!}{ 
    \begin{tabular}{l|c|l|l|c|l}     
    \mytoprule
    Distribution    & $l$-MDA & \multicolumn{2}{c|}{Normalization constants} & $p$-MDA & Normalization constants \\\hline
       Log-Fr{\'e}chet& - & \multicolumn{2}{c|}{-}&$ H_{1,\alpha}  $ &\tabincell{l}{$\alpha_{n}=1$ \\ $\beta_{n}=n^{-1/\alpha}$} \\\hline
       Log-Polynomial& $\Psi_{\alpha } $ &\multicolumn{2}{l|}{\tabincell{l}{$a_n=n^{1/\alpha}/\expon{\gamma(F)}$\\$ b_n=\expon{\gamma(F)}$}}  &$H_{2,\alpha}  $ &\tabincell{l}{$\alpha_{n} =\expon{-n^{1/\alpha}\gamma(F)} $\\$ \beta_{n} = n^{1/\alpha}$ }  \\\hline
        $\text{U}[l,u]$& $\Psi_1$ &\multicolumn{2}{l|}{\tabincell{l}{$a_n=\frac{n}{u-l}$\\$ b_n=u $} } &$H_{2,1}  $ &\tabincell{l}{$\alpha_{n} =u^{-un/(u-l)} $\\$ \beta_{n} = \frac{un}{u-l}$ }  \\\hline
       standard normal & $\Lambda$ & \multicolumn{2}{l|}{\tabincell{l}{$a_n=(2\log n)^{1/2} $\\$ b_n= a_n-\frac{\displaystyle\log ( 4 \pi\log n)}{\displaystyle 2 a_n}$}}
       & $H_5$ & \tabincell{l}{$\alpha_n=b_{n}^{-a_{n}b_{n}} $\\$ \beta_n=a_{n}b_{n}$} \\\hline
              General error & $\Lambda$ &\multicolumn{2}{c|}{\tabincell{l}{$a_n= 2^{-1/\nu}(\nu/\lambda)(\log n)^{1-1/\nu}$ \\ $b_n = \frac{\displaystyle\nu\log n -\frac{\nu-1}\nu \log \left[2\Gamma\fracl1\nu\log n\right]}{\displaystyle a_n}$}}  & $H_5$ &\tabincell{l}{$\alpha_n = b_{n}^{-a_{n}b_{n}} $\\$ \beta_n=a_{n}b_{n}$} \\\hline
       Fr{\'e}chet/Pareto & $\Phi_\alpha$ &\multicolumn{2}{l|}{\tabincell{l}{$a_n=n^{-1/\alpha}$ \\ $b_n=0$}} &$H_5$ & \tabincell{l}{$\alpha_n=1/n
       $\\$ \beta_n= \alpha 
       $}\\ 
       \hline
       \multirow{2}{*}[-1ex]{Skew-Normal} & \multirow{2}{*}[-1ex]{$\Lambda$} & $\lambda \geq 0$ & \tabincell{l}{$a_n=b_n$\\$ b_n= F_{\lambda}^{-1}(1-1/n)$ } & \multirow{2}{*}[-1ex]{$H_5$} & \tabincell{l}{$\alpha_n=a_n^{-a_n^{2}}$ \\$ \beta_n=a_n^{2}$}\\ \cline{3-4} \cline{6-6}
       & & $\lambda<0$ & \tabincell{l}{$a_n=(1+\lambda^2)b_n$\\$ b_n=F_{\lambda}^{-1}(1-1/n)$}& & \tabincell{l}{$\alpha_n=b_n^{-a_n b_n}$ \\$ \beta_n=a_n b_n$} \\\hline
       \multirow{3}{*}[-3ex]{Polynomial} & \multirow{3}{*}[-3ex]{$\Psi_\alpha$} & $\gamma(F)>0$ &\multirow{3}{*}[-3ex]{\tabincell{l}{$a_n=n^{1/\alpha}$\\ $b_n=\gamma(F)$}} &  $H_{2,\alpha}$ & \tabincell{l}{$\alpha_n=\gamma(F)^{-\gamma(F) n^{1/\alpha}} $\\$ \beta_n=\gamma(F) n^{1/\alpha}$}\\ \cline{3-3} \cline{5-6}
       & & $\gamma(F)=0$ & & $H_6$ & \tabincell{l}{$\alpha_n= n^{1/\alpha}$\\ $\beta_n=1/\alpha$} \\ \cline{3-3} \cline{5-6}
       & & $\gamma(F)<0$ & &  $H_{4,\alpha}$ & \tabincell{l}{$\alpha_n=(-\gamma(F))^{\gamma(F) n^{1/\alpha}} $\\$ \beta_n=-\gamma(F) n^{1/\alpha}$}\\
       \mybottomrule
    \end{tabular}}
       \label{Table: CommonDistribution}
\end{table}

\section{Numerical Studies} 
\label{sec: simulation}
We conduct Monte Carlo simulations to illustrate Theorem \ref{TheoremMaxMax} with examples provided in the previous section. Throughout the simulation, repeated time  $m = 10^4$ unless otherwise stated.  Section \ref{sec:datafitting} illustrates the efficient approximations of  theoretical distributions.  
Section \ref{sec: convergence} compares the convergence rates of extremes of competing risks under power and linear normalization. Finally, Section \ref{Competing Detection} will demonstrate the method to detect competing risks, which will be applied to real data in Section \ref{sec:Real data}. 
\subsection{Efficiency of Distribution Approximations}\label{sec:datafitting}
We will illustrate the efficiency of  distribution approximation, using a graphical procedure of the deviation between theoretical density curves and the histograms of simulated samples. 

\textbf{Case I: Convergence to accelerated $p$-max stable distributions.} 
The accelerated case of limit distribution for maxima of maxima, i.e.,  $H^{(1)}(A|x|^B \sign(x))H^{(2)}(x)$ introduced in Theorem \ref{TheoremMaxMax}(i),
usually occurs when the limit distributions of concordant sub-maxima are of the same type. Meanwhile, the sub-maxima with lighter-tailed limit distribution should have a larger sample size, so that strong competing risk prevents one sub-maxima from dominating the other. Typical examples are log-Fr{\'e}chet/Pareto \& log-Fr{\'e}chet/Pareto, log polynomial \& log polynomial, and Pareto \& Pareto given in Examples \ref{FandF}, \ref{WandW} and \ref{ex3.8}, respectively.

\begin{figure}[h]
        \centering
\subfigure[]{
\includegraphics[width=0.45\linewidth]{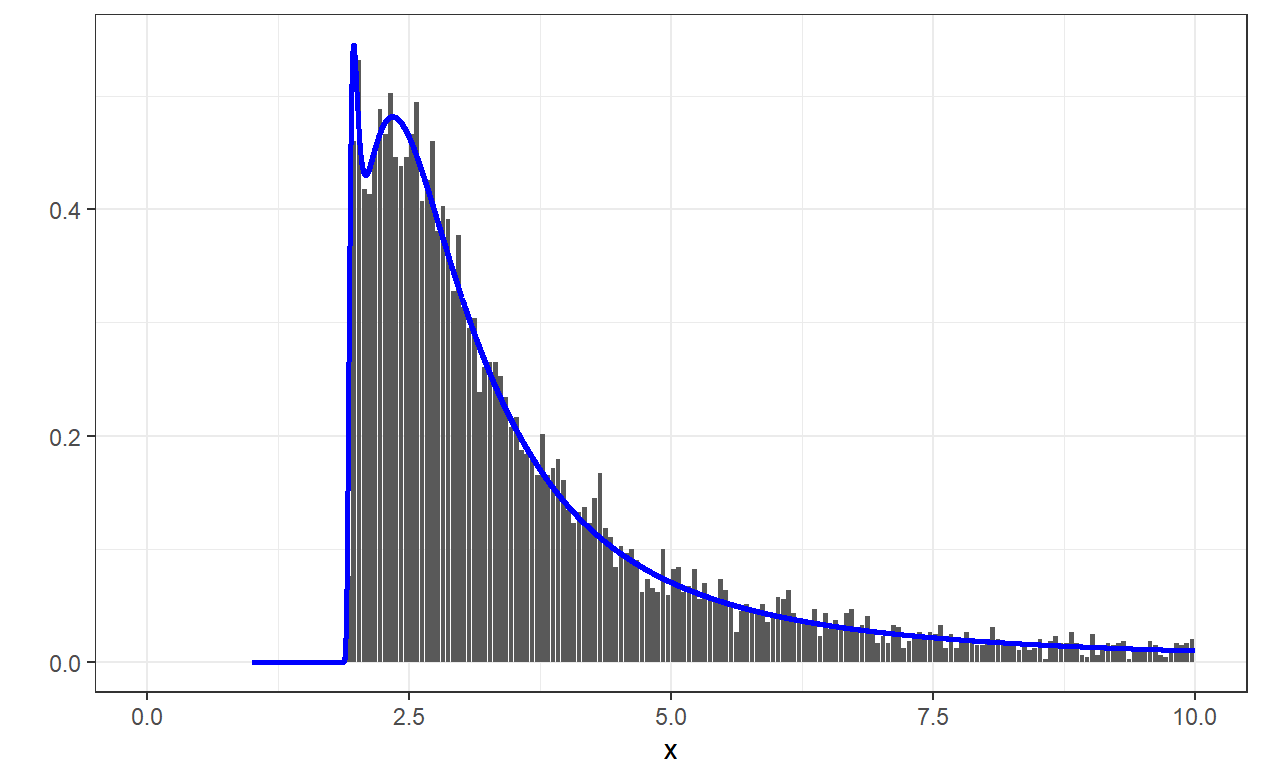}
}
\subfigure[]{
\includegraphics[width=0.45\linewidth]{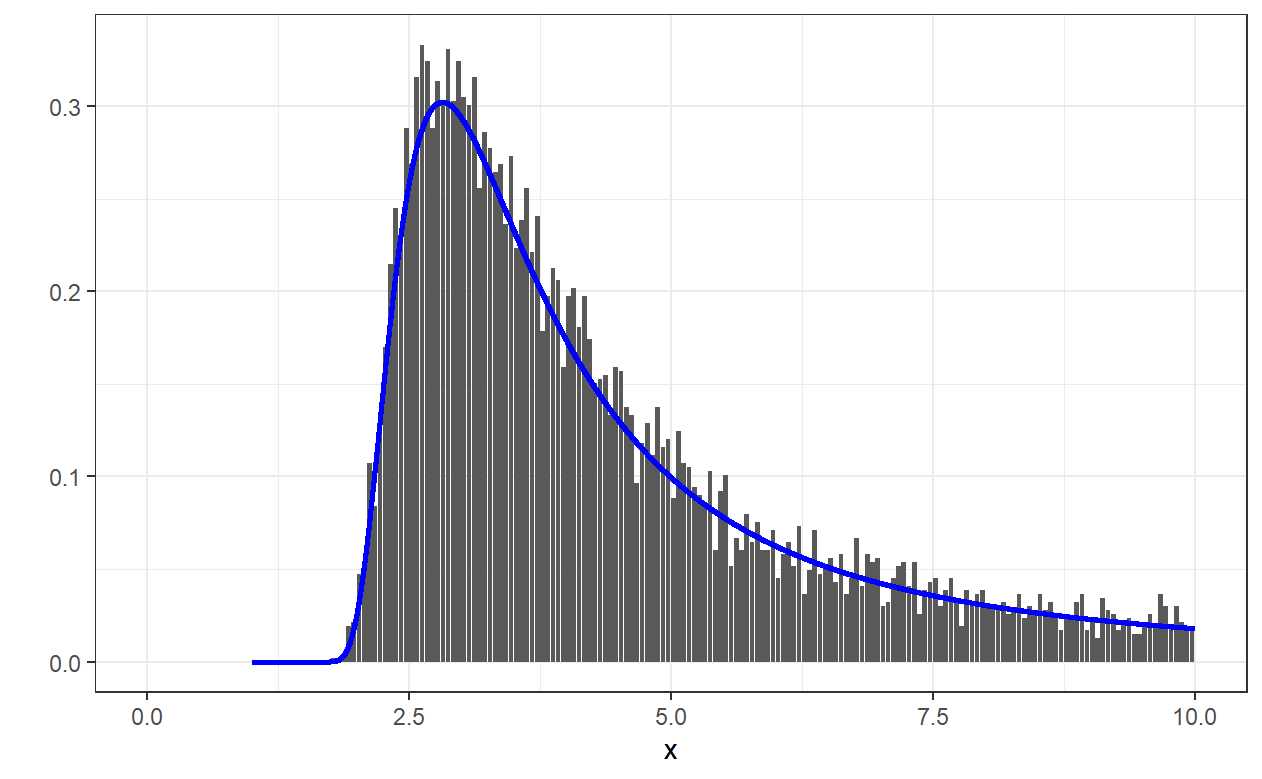}
}
\subfigure[]{
\includegraphics[width=0.45\linewidth]{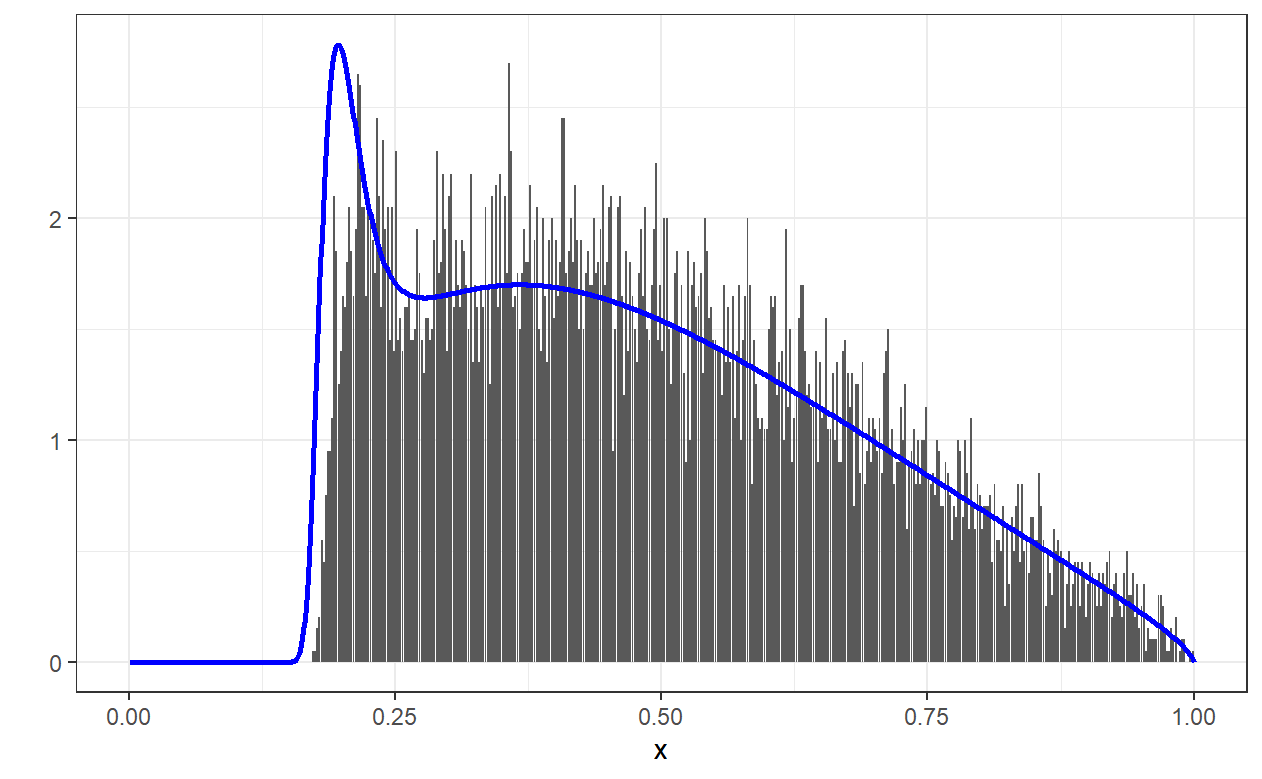}
}
\subfigure[]{
\includegraphics[width=0.45\linewidth]{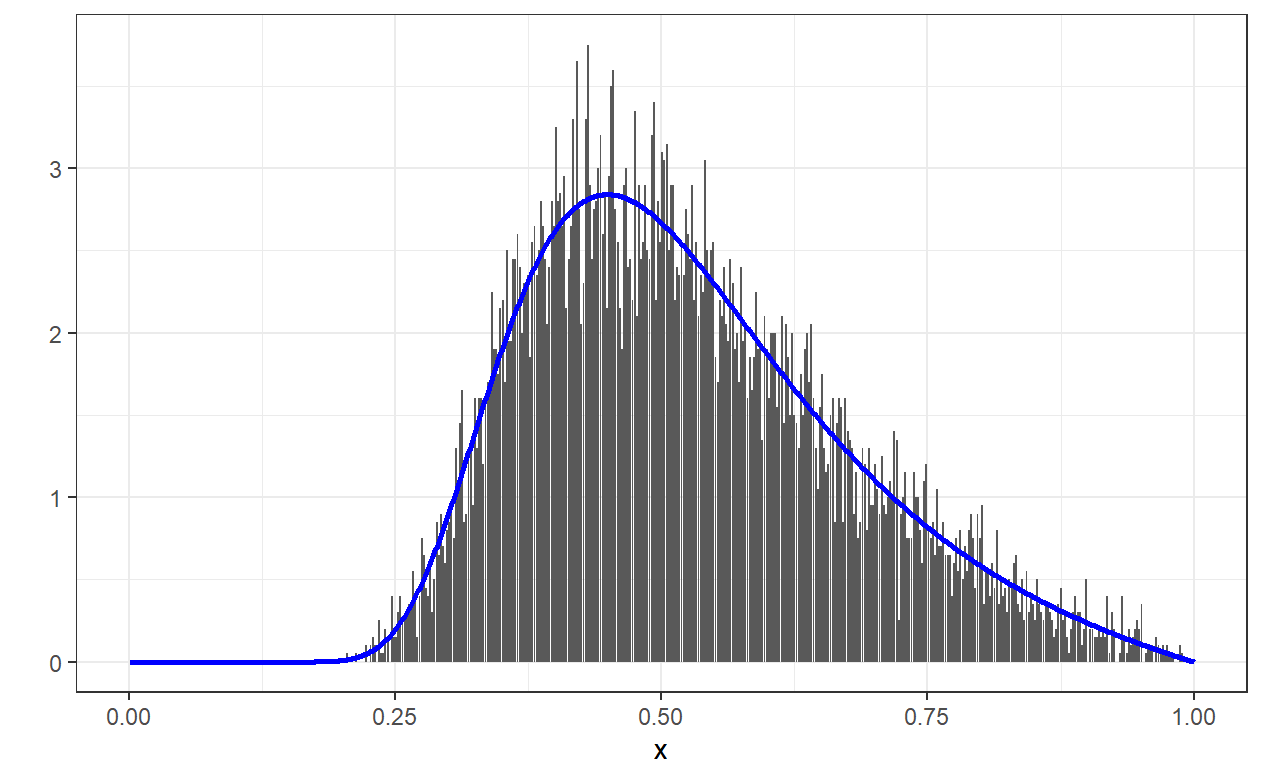}
}
    \caption{Histogram of power normalized $M_n = \max(M_{1, n_1}, M_{2, n_2})$ with both $M_{j,n_j}$'s from log-Fr\'echet and log-polynomial. The parameters and sample size are taken as $(\alpha_1,\alpha_2,n_1,n_2)$ equal $(40,3,10^5,8)$ in (a), and $ (4,2,10^4,100)$ in (b,d), and $(20,1.7,10^5,6)$ in (c). 
    The curves denote the densities of  $H_{1, \alpha_1}(|x|^B \sign(x))H_{1, \alpha_2}(x)$ in (a,b) with $B=1.5,1$, respectively and $H_{2, \alpha_1}(|x|^B \sign(x))H_{2, \alpha_2}(x)$ in (c,d) with with $B=0.6,1$, respectively.}  
    \label{Figacctype}
\end{figure}

Figure \ref{Figacctype} \cl{represents} 
the accelerated cases with samples generating according to Examples \ref{FandF} and \ref{WandW} respectively for (a,b) and  (c,d).  We take  $(\alpha_1, \alpha_2) = (40,3), (n_1, n_2)=(10^5, 8)$, and $(\alpha_1,\alpha_2) = (4,2), (n_1, n_2)=(10^4,100)$ in Figure \ref{Figacctype}(a,b), respectively. Both sample sizes satisfy that $n_2 =B^{\alpha_2}n_1^{\alpha_2/\alpha_1}$ with $B=1.5, 1$ for Figure \ref{Figacctype}(a,b), respectively. According to Theorem \ref{TheoremMaxMax}(i), the approximated distribution is given by
\begin{eqnarray*}
   H(x) &=&H_{1, \alpha_1}(|x|^B \sign(x))H_{1, \alpha_2}(x) \\&=&\expon{ - [(\log x)^{-\alpha_2} + (B\log x)^{-\alpha_1}]},\quad x>1.
\end{eqnarray*}
Clearly, we see the nice approximation of the density of $H(x)$ to the histogram of the sample max of max in Figure \ref{Figacctype}(a,b).  While in Figure \ref{Figacctype}(c,d), we generate samples from log-polynomial \& log-polynomial  according to Example \ref{WandW}, the parameters and sample sizes are $(\alpha_1,\alpha_2,n_1,n_2)= (20,1.7,10^5,6)$ in (c), and $
(4,2,10^4,100)$ in (d). Both sample sizes satisfy $n_2 =B^{-\alpha_2}n_1^{\alpha_2/\alpha_1}$ with $B=0.6, 1$, respectively. According to Theorem \ref{TheoremMaxMax}(i), the approximated distribution is given by
\begin{eqnarray*}
    H(x) &=& H_{2, \alpha_1}(|x|^B \sign(x))H_{2, \alpha_2}(x) \\&=& \expon{ - [(-\log x)^{\alpha_2} + (-B\log x)^{\alpha_1}]},\quad 0<x<1.
\end{eqnarray*}
Note that the result shows the distribution convergence of the competing risks to both bimodal (a,c) and unimodal (b,d) accelerated $p$-max stable law, respectively. While all of the four distributions cannot be reduced to any one of the six $p$-max stable distribution types. In addition, better goodness of fit was given for log-Fr{\'e}chet \& log-Fr{\'e}chet competing risks in comparison with that  with log-polynominal \& log-polynomial case. This is likely due to the fact that the log-Fr{\'e}chet is an exact $p$-max stable distribution, resulting in faster convergence of the simulation process.

\textbf{Case II. Convergence to $p$-max stable distributions}. According to Theorem \ref{TheoremMaxMax}(ii), apart from some special cases of competing risks, in most circumstances, one sub-maxima $(M_{2,n_2})$ can be  greater than the other sub-maxima $(M_{1, n_1})$, leading to the maxima of maxima $(M_n)$ being exactly the greater sub-maxima $(M_{2,n_2})$ almost surely. We call this case as $M_{2,n_2}$ dominates $M_{1, n_1}$, moreover, we say that the distribution $F_2(x)$ (usually heavier-tailed) dominates the distribution $F_1(x)$ (usually lighter-tailed). Note that the dominant situation may be weakened or even reversed with respect to different sample sizes. However, the dominance can always exist when the lighter-tailed distribution has a finite right-endpoint. For example, $U\left [ 2,4 \right ]$ with a smaller right-endpoint is dominated by $U\left [ 1,5 \right ]$. Similarly, the polynomial growth distribution ($\alpha_1=4$) with finite right-endpoint is dominated by the Fr\'echet distribution ($\alpha_2=2$).

\begin{figure}[h]
    \centering
\subfigure[]{
\includegraphics[width=0.45\linewidth]{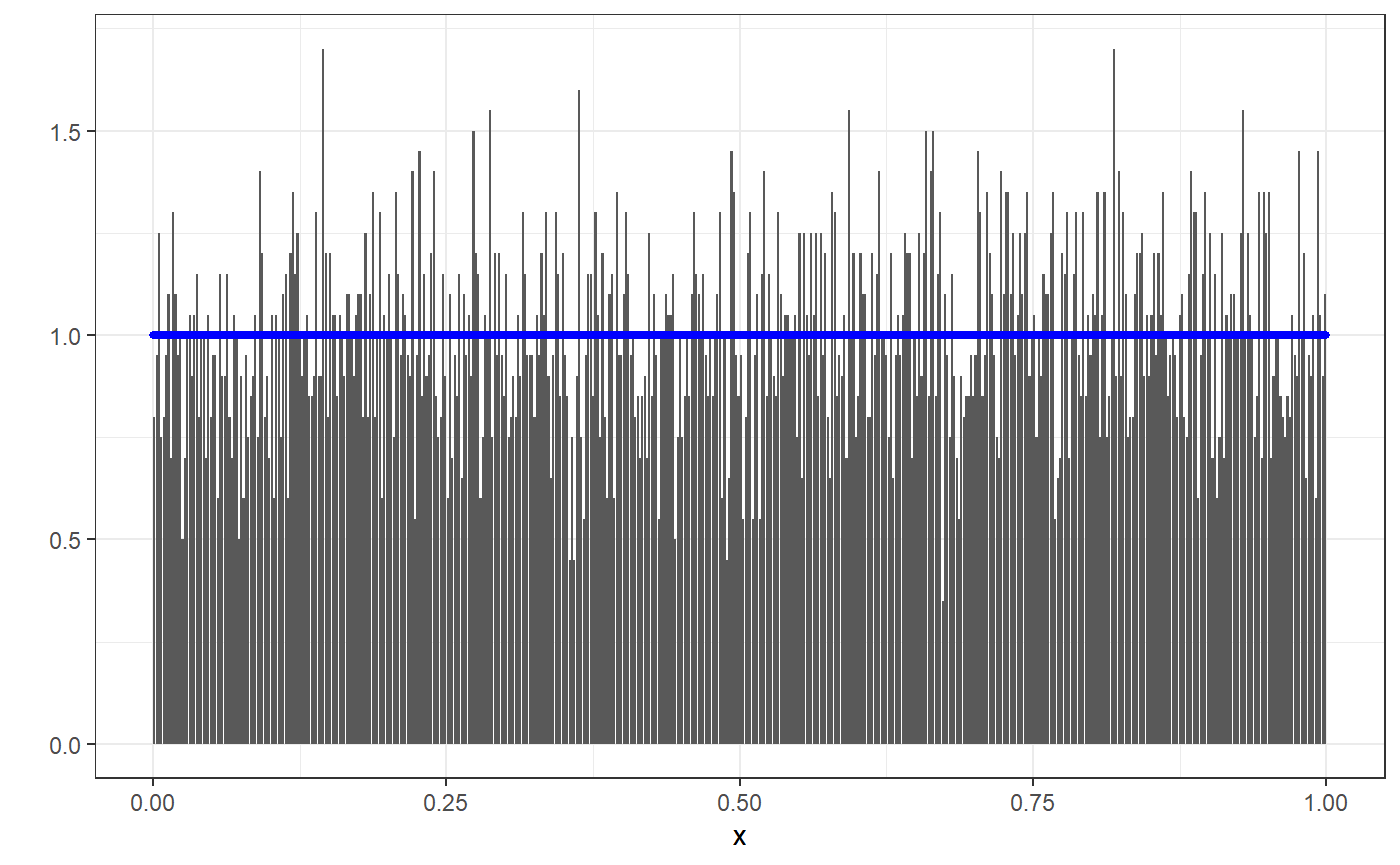}
}
\subfigure[]{
\includegraphics[width=0.45\linewidth]{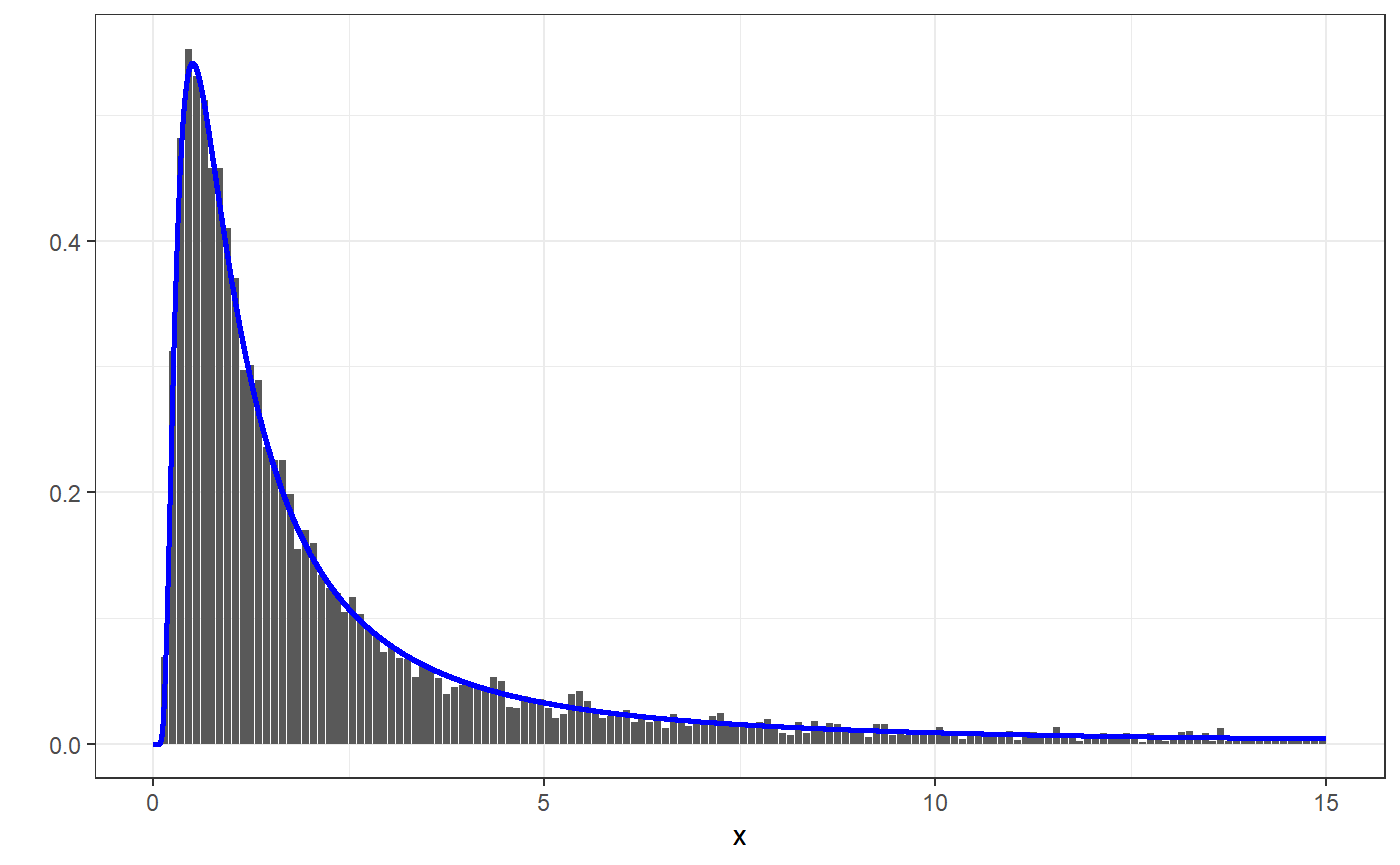}
}

\caption{Histogram of power normalized $M_n = \max(M_{1, n_1}, M_{2, n_2})$
with  $M_{j,n_j}$'s from $U[2,4]$ \& $U[1,5]$ in (a) and polynomial growth \& Fr\'echet ($(\alpha_1, \alpha_2)=(4,2)$) in (b). Both sample sizes are $n_1=n_2=100$. The density curves correspond to $H_{2,1}$ (standard uniform)  and $H_5$ (unit Fr\'echet) in (a,b), respectively.}
    \label{Figdominatetype}
\end{figure}

Figure \ref{Figdominatetype}(a,b) shows distribution approximation for competing risks from $U\left [ 2,4 \right ]$ \& $U\left [ 1,5 \right ]$, and polynomial growth ($\alpha_1=4$) \& Fr\'echet ($\alpha_2=2$), with same sample size of $n_1 = n_2 = 100$. Note in both examples, $F_2(x)$ dominates  $F_1(x)$, and the limit distributions of $M_n$ become $H(x)=H^{(2)}(x)=H_{2,1}(x)$ in (a) and $H(x)=H^{(2)}(x)=H_5(x)$ in (b), which are exactly standard uniform distribution and unit Fr\'echet distribution, respectively.

\textbf{Case III. Convergence to left-truncated $p$-max stable distribution}. The case of left-truncated $p$-max stable distribution is usually due to the limit distributions of sub-maxima in different types of $p$-MDA but with overlapping supports, see Theorem \ref{remark1}(iii). In this case, the lighter-tailed limit distribution degenerates at an interior point of the heavier-tail distribution, resulting in a left-truncated limit distribution of $M_n$. As mentioned in the examples of Pareto \& log-Fr{\'e}chet (Example \ref{PandF}) and general error \&  log-Fr{\'e}chet (Example \ref{GandF}), the limit distributions of $M_n$ become $H^{(2)}(x)\I{x>x_0}$ with the jump point $(x_0, H^{(2)}(x_0))$, where $H^{(2)}(x)$ is the heavier-tailed distribution. The corresponding likelihood becomes $H^{(2)}(x_0)\I{x=x_0}+h^{(2)}(x)\I{x>x_0}$, where 
$h^{(2)}(x)$ denotes the density of $H^{(2)}(x)$. 

\begin{figure}[htbp!]
    \centering
\subfigure[]{
\includegraphics[width=0.45\linewidth]{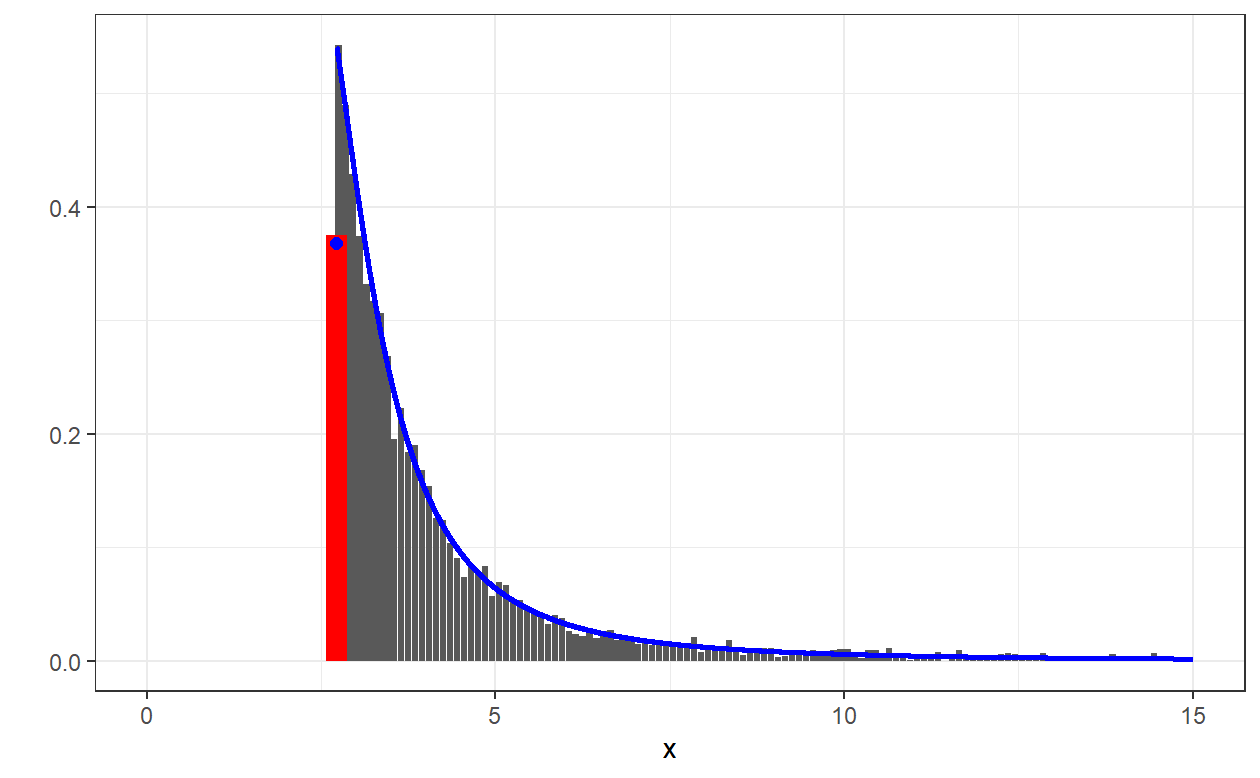}
}
\subfigure[]{
\includegraphics[width=0.45\linewidth]{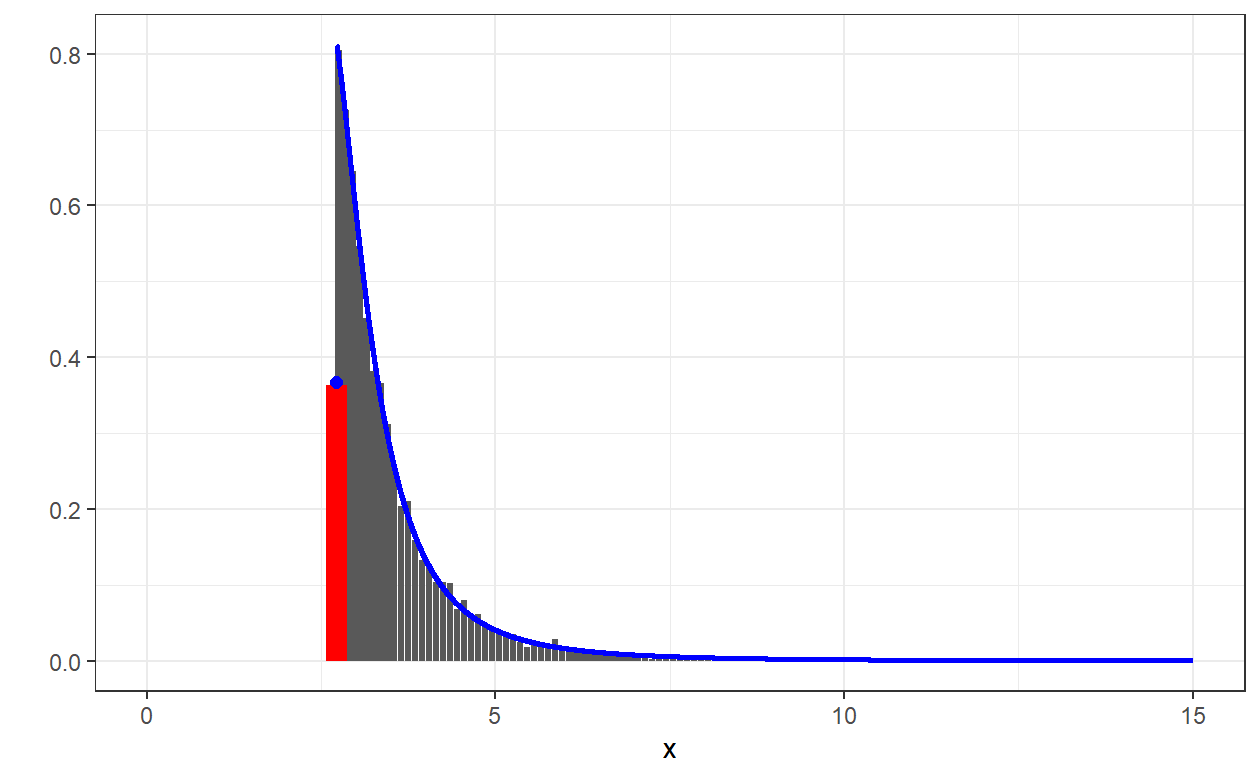}
}

\caption{Histogram of power normalized $M_n = \max(M_{1, n_1}, M_{2, n_2})$ from Pareto \& \kai{log-Fr\'echet} with $(\alpha_1,\alpha_2,n_1,n_2)=(2,4,10^4,449)$ in (a) and General error \& \kai{log-Fr\'echet} with $(\nu,\alpha_2,n_1,n_2)=(1,6,10^4,120)$ in (b).
The likelihood curves are $H_{1,\alpha_2}(x_0)\I{x=x_0}+h_{1,\alpha_2}(x)\I{x> x_0}$ with $x_0 =\e$ and $\alpha_2 =4,6$, respectively. Both  jump points  are $(\mathrm{e}, \mathrm{e}^{-1})$ with empirical step-length  of  0.363 and 0.365 in (a,b), respectively.}

\label{Figmixtype}
\end{figure}

We conduct simulation for Pareto \& log-Fr\'echet example with parameters and sample sizes $(\alpha_1,\alpha_2,n_1,n_2)=(2,4,10^4,\hu{449})$ in Figure \ref{Figmixtype}(a), and for general error \& log-Fr\'echet with $(\nu,\alpha_2,n_1,n_2)=(1,6,10^4,120)$  in \ref{Figmixtype}(b). Thus, the left-truncated $p$-max stable distribution type limit distribution of $M_n$ are $H_{1,\alpha_2}(x)\I{x> x_0}$ with corresponding likelihood  $H_{1,\alpha_2}(x_0)\I{x=x_0}+h_{1,\alpha_2}(x)\I{x> x_0}$ with $x_0 =\e, \e^{1/\nu}$ and $\alpha_2 = 4,6$ in (a,b), respectively.

The isolated points in both Figure \ref{Figmixtype}(a,b) are $(\mathrm{e}, H_{1,\alpha_2}(\mathrm{e})) = (\mathrm{e}, \mathrm{e}^{-1})$, given by the theoretical probabilities of jump points with height $1/\e=0.368$, where the red bars denote the sample mean of $\mathbb I\{ \alpha_{2,n_2}|M_n|^{\beta_{2,n_2}}\sign(M_n) \leq \mathrm{e}\}$ equals to 0.363 and 0.365 in (a,b), respectively.

\subsection{Comparison of Convergence Rates} \label{sec: convergence}

In this section, we compare the convergence rate of maxima of maxima under linear normalization and power normalization for those distributions belonging to both  $l$-MDA and $p$-MDA, see e.g., Examples \ref{WandW}, \ref{UandU}, \ref{ex3.8} and \ref{ex3.9}. To this end, we apply the three goodness-of-fit tests with implementation in R-stats package \citep{Rstats} and R-goftest package \citep{Rgoftest}. 

In Table \ref{Table2}, we generate $m=10^4$ samples of size $n_1=n_2=10^4$ for each pair of competing risks. The reported test statistics and its $p$-values show that, a better convergence rate of power normalized extreme is given by Example \ref{WandW} (log-polynomial \& log-polynomial), while similar convergence speed of linear and power normalization in Examples \ref{ex3.8} (Pareto \& Pareto) and \ref{ex3.9} (polynomial \& Fr{\'e}chet).

\begin{table}[H]

    \centering
      \caption{Comparison of convergence rate of the specified competing risks under power normalization and linear normalization. Here both sample size  $n_1=n_2=10^4$.}
      \label{Table2}
      \setlength{\tabcolsep}{1mm}
    \begin{tabular}{c|c|c|c|c|c}
    \mytoprule
     Distributions & Parameter ($\alpha$) & Normalization  & KS test & CVM test & AD test\\\hline
     log-Pareto &  4 & $l$-max & $0.0095 (0.33)$ & $0.1649 (0.35)$ & $0.8164 (0.47)$\\\cline{2-6}
     log-Pareto &  2               & $p$-max & $\mathbf{0.0075 (0.63)}$ & $\mathbf{0.0934 (0.62)}$ & $\mathbf{0.5451 (0.70)}$ \\\hline                  
    Pareto & 4& $l$-max & $0.0075 (0.63)$ & $0.0934 (0.62)$ & $0.5451 (0.70)$\\\cline{2-6} 
      Pareto& 2   & $p$-max & $0.0075 (0.63)$ & $0.0934 (0.62)$ & $0.5451 (0.70)$ \\\hline
     Polynomial& 4 & $l$-max & $0.0067 (0.77)$ & $0.0640 (0.79)$ & $0.3881 (0.86)$\\\cline{2-6} 
      Fr{\'e}chet & 2& $p$-max & $0.0067 (0.77)$ & $0.0640 (0.79)$ & $0.3881 (0.86)$\\
    \mybottomrule
    \end{tabular}

\end{table}
In Table \ref{Table3}, the data is generated according to Example \ref{UandU} ($U [2, 4]$ \& $U [1, 5]$) with sample size satisfying $n_2= n_1 + 50$ and $n_1 = 100, 200, 300$. We see that the model-based approximation is confirmed with stronger evidence with generated samples of larger sample size. Meanwhile,  the approximation of power normalized extremes performs always better than linear normalization.

\begin{table}[H]
 \centering
      \caption{Comparison of convergence rate of maxima of maxima from competing risks of $U[2,4]\&U[1,5]$ under power and linear normalization. Here sample size satisfies $n_2= n_1+50$ and repeated $m=10^4$.}
      \setlength{\tabcolsep}{3mm}
    \begin{tabular}{c|c|c|c|c}
    \mytoprule
    Sample size ($n_1$) & Normalization  & KS test & CVM test & AD test\\\hline
     \multirow{2}{*}{$100$} & $l$-max & $0.0130 (0.07)$ & $0.3053 (0.13)$ & $1.6542 (0.14)$\\\cline{2-5}
                        & $p$-max & $\mathbf{0.0119 (0.12)}$ & $\mathbf{0.2396 (0.20)}$ & $\mathbf{1.2699 (0.24)}$
                        \\\hline                  
     \multirow{2}{*}{$200$} & $l$-max & $0.0100 (0.27)$ & $0.2837 (0.15)$ & $1.9364 (0.10)$\\\cline{2-5} 
                         & $p$-max & $\mathbf{0.0091 (0.38)}$ & $\mathbf{0.2418 (0.20)}$ & $\mathbf{1.6636 (0.14)}$ \\\hline
     \multirow{2}{*}{$300$} & $l$-max & $0.0073 (0.67)$ & $0.1071 (0.55)$ & $0.8259 (0.46)$\\\cline{2-5} 
                         & $p$-max & $\mathbf{0.0072 (0.67)}$ & $\mathbf{0.0954 (0.61)}$ & $\mathbf{0.7574 (0.51)}$\\
    \mybottomrule
    \end{tabular}
    \label{Table3}
\end{table}

\subsection{Detection of Competing Risks}\label{Competing Detection}
Competing risks and heterogeneity may exist among extreme observations that are originated from competing sources. This section aims to detect competing risks by comparing the performance of accelerated and single $p$-max stable models. Firstly, we generate two samples, one represents the dominant case and the other represents the competing case. Subsequently, we obtain the maximum likelihood estimators for parameters involved, and then apply goodness-of-fit tests and likelihood ratio test (LRT). If both theoretical distributions fit the sample comparably well, we may conclude no competing risks. However, if the accelerated type performs significantly better than the single case, heterogeneity probably exists.

We will generate a sample of size $m=10^4$ from accelerated $p$-max stable model. In particular, we have
$$Z_i = \max(Z_{i1}, Z_{i2}), \quad i=1, \ldots, m$$
with both $Z_{i1}$ and $Z_{i2}$ following $p$-max stable d.f. $H_1^{\xi}(z; \mu,  \sigma)$ with specified values of $\boldsymbol\theta = (\mu, \sigma,\xi)$. 
To illustrate the dominated case, we set both location ($\mu$) and  shape ($\xi$) parameters larger in the generating distribution for one sub-sample so that it can dominate the other sub-sample. In particular, we generate sub-samples from $p$-max stable distributions with parameters $(\mu_1,\sigma_1,\xi_1)=(2,1,-0.2)$ and $(\mu_2,\sigma_2,\xi_2)=(0,1,-1)$.
While for the competing condition, we associate the larger location parameter with smaller shape parameter, with specific choices of $(\mu_1,\sigma_1,\xi_1)=(3,1,0.1)$ and $(\mu_2,\sigma_2,\xi_2)=(2,1,0.5)$ to generate sub-samples.

Table \ref{Table4} provides the maximum likelihood estimates and goodness-of-fit for both dominated and competing cases. In the dominated case, although both models are well-fitted and make it difficult to confirm the competing risks, three goodness-of-fit tests and the LRT test ($p<0.001$) still illustrate the significant improvement of accelerated $p$-max stable model. However, in the competing case, the single $p$-max stable model failed in all three goodness-of-fit tests at the 5\% level, while the accelerated $p$-max stable model shows remarkably better fitness, suggesting the existence of strong competing risks (heterogeneity). 

\begin{table}[H]
    \centering
    \caption{Maximum likelihood estimates of parameters (standard error), goodness-of-fit test and LRT statistic ($p$-value) based on \cl{accelerated $p$-max stable} and single $p$-max stable models. Here the data is a random sample of size $m=10^4$ from \cl{$H_1^{\xi_1}(x;\mu_1, \sigma_1)H_1^{\xi_2}(x;\mu_2, \sigma_2)$} with specified parameters. 
    }
     \resizebox{\linewidth}{!}{ 
    \begin{tabular}{c|r|r|r|r} \mytoprule
    \multirow{3}{*}{-}& \multicolumn{2}{c|}{$(\mu_1,\sigma_1,\xi_1)=(2,1,-0.2)$}& \multicolumn{2}{c}{$(\mu_1,\sigma_1,\xi_1)=(3,1,0.1)$} \\
    &\multicolumn{2}{c|}{$(\mu_2,\sigma_2,\xi_2)=(0,1,-1) \ $}& \multicolumn{2}{c}{$(\mu_2,\sigma_2,\xi_2)=(2,1,0.5)$} \\ \cline{2-5}
    & accelerated $p$-max type & $p$-max stable distribution & accelerated $p$-max type & $p$-max stable distribution\\ \hline
    $\widehat{\mu}_1$  & $2.00(0.01)$ & $2.00(0.01)$ & $2.80(0.15)$ & $3.38(0.01)$ \\\hline
    $\widehat{\sigma}_1$ & $1.00(0.01)$ & $0.95(0.01)$  & $0.97(0.03)$  & $1.24(0.01)$ \\\hline
    $\widehat{\xi}_1$ & $-0.20(0.01)$  & $-0.17(0.01)$  & $0.05(0.05)$  & $0.31(0.01)$\\\hline
    $\widehat{\mu}_2$ & $0.29(0.01)$ & -  & $2.41(0.25)$   & -\\\hline
    $\widehat{\sigma}_2$ & $0.58(<0.005)$ & - & $1.07(0.07)$   & -\\\hline
    $\widehat{\xi}_2$ & $-0.81(0.01)$ & - & $0.47(0.03)$ & -\\\hline
    KS & $\mathbf{0.0078(0.58)}$ & $0.0114(0.15)$ & $\mathbf{0.0062(0.82)}$ & $0.0200(<0.005)$ \\\hline
    CVM & $\mathbf{0.0824(0.62)}$ & $0.2133(0.24)$ & $\mathbf{0.0337(0.96)}$ & $1.0992(<0.005)$ \\\hline
    AD & $\mathbf{0.4923(0.75)}$ & $2.1149(0.08)$ & $\mathbf{0.2565(0.97)}$ & $7.0569(<0.005)$ \\\hline
    LRT & \multicolumn{2}{c|}{$64.9010(<0.001)$}& \multicolumn{2}{c}{$119.2000(<0.001)$}\\
    \mybottomrule
    \end{tabular}
    }
    \label{Table4}
\end{table}

\section{Real Data Analysis}
\label{sec:Real data}

In real data analysis, the accelerated $p$-max/$p$-min stable model can be applied even when the single $p$-max/$p$-min stable model is not applicable because of the substantial heterogeneity and strong competing risks. For example, the survival time of patients who get heart transplanted may be affected by complications and comorbidities, which can be regarded as competing risks leading to the infeasibility of single $p$-min stable model. With the atmospheric (ground-level ozone) and survival time (Stanford heart plant) datasets, we utilize both goodness-of-fit tests and diagnostic plots to show the robustness of accelerated $p$-max/$p$-min stable models, as well as investigate the potential heterogeneity by model comparison through LRT.

\begin{figure}[H]
        \centering
\subfigure[]{
\includegraphics[width=0.4\linewidth]{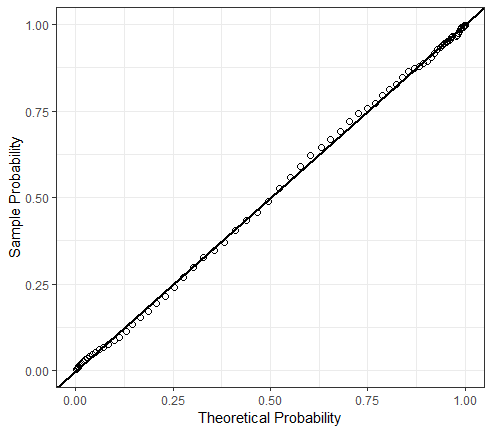}
}
\subfigure[]{
\includegraphics[width=0.4\linewidth]{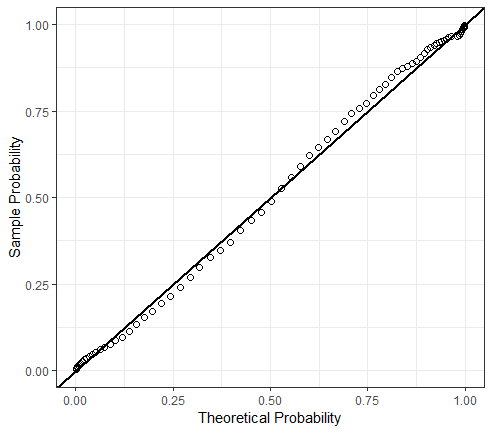}
}
\subfigure[]{
\includegraphics[width=0.4\linewidth]{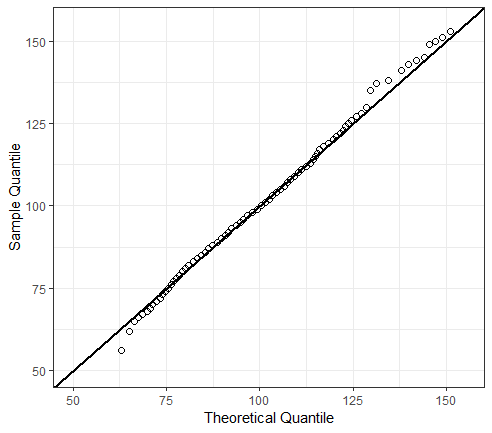}
}
\subfigure[]{
\includegraphics[width=0.4\linewidth]{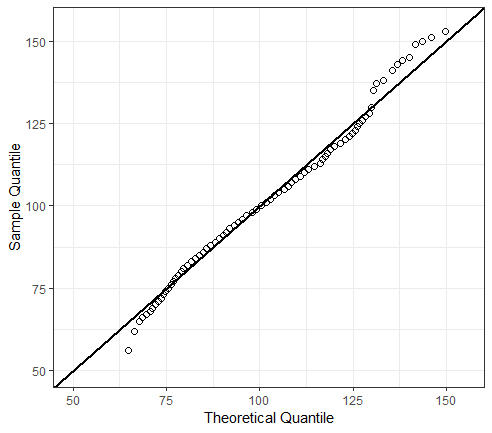}
}
    \caption{P-P plots and Q-Q plots  for Ground-level ozone data \cl{based on} accelerated $p$-max stable models in (a,c) and single $p$-max stable models in (b,d). Relevant parameters are  given in Table \ref{Table5}.} 
    
    \label{Diagostic plot1}
\end{figure}

The ground-level ozone dataset \citep{extRemes} contains the maximum to the fourth-highest ozone readings in 1997 from 513 monitoring stations in the eastern United States, see details in \citet{Fuentes2003}. As for the annual maxima, compared to the single $p$-max stable model, the three goodness-of-fit tests in Table \ref{Table5} report better overall fitness (with larger $p$-values) for the accelerated $p$-max stable model. The accelerated model also exhibits especially superior tail fitness (especially the right tail) by the P-P plot and Q-Q plot (Figure \ref{Diagostic plot1}). This conclusion is further supported by the likelihood ratio test $(p<0.001)$, and all the above findings collectively suggest the existence of competing risks (heterogeneity) among the maximum observations from different stations.

\begin{figure}[H]
        \centering
\subfigure[]{
\includegraphics[width=0.4\linewidth]{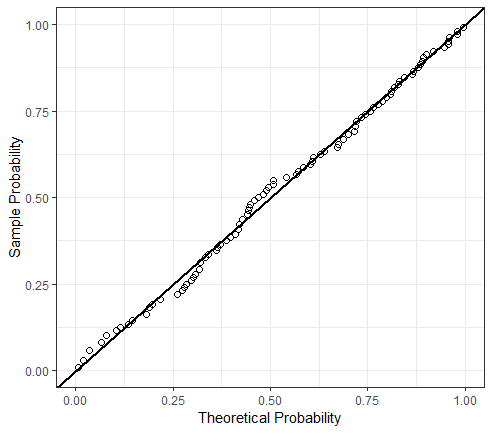}
}
\subfigure[]{
\includegraphics[width=0.4\linewidth]{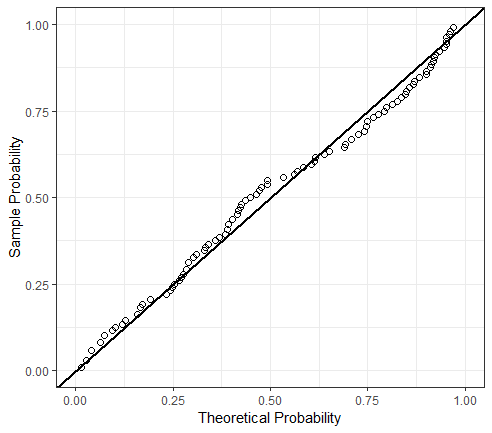}
}
\subfigure[]{
\includegraphics[width=0.4\linewidth]{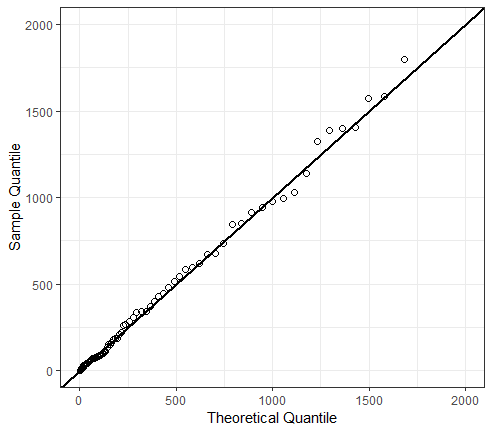}
}
\subfigure[]{
\includegraphics[width=0.4\linewidth]{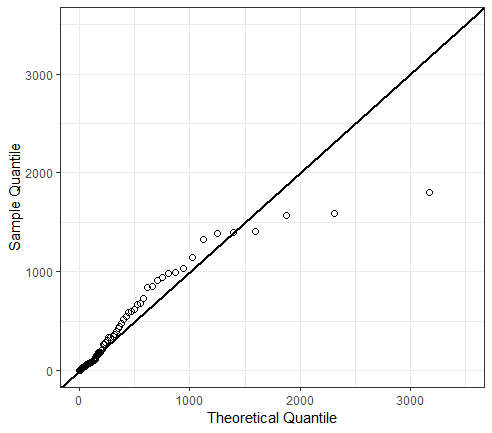}
}
    \caption{P-P plots and Q-Q plots  for Stanford heart plant data \cl{based on} accelerated $p$-min stable models in (a,c) and single $p$-min stable models in (b,d). Relevant parameters are  given in Table \ref{Table5}.
    } 
    
    \label{Diagostic plot2}
\end{figure}

The second dataset, Stanford heart transplant dataset \citep{Stanford} provides information about the survival time of 103 patients from their acceptance into the transplant program until 1st April 1974. Among the 103 patients, 69 received heart transplants, 75 patients died, and the remaining 28 patients were censored. Since the survival times of the patients are possibly affected by commodities and complications, in this scenario, we apply the accelerate $p$-min stable model and single $p$-min stable model. Upon analyzing the survival times of all 103 patients, despite good fitness indicated by the goodness-of-fit tests (Table \ref{Table5}) for the single $p$-min stable model ($p > 0.80$), the accelerate-type model still outperforms ($p > 0.95$) and the significant improvement is confirmed by LRT $(p=0.005)$. More importantly, taking the effects of commodities and complications into account, the accelerate-type model provides shorter estimated survival time and exhibits superior right-tail fitness, as demonstrated by the P-P plot and Q-Q plot in Figure \ref{Diagostic plot2}. Based on these findings, it is possible to obtain some evidence about the existence of relatively weak competing risks from milder commodities and complications.

The standard errors for maximum likelihood estimators in the accelerated cases are usually larger than the ones in the single $p$-max or $p$-min case (Table \ref{Table5}), because the accelerated cases have three more parameters to estimate. Hence this approach is more appropriate and promising to detect the competing risk for a large amount of clustering data.

\begin{table}[H]
    \centering
    \caption{Maximum likelihood estimates of parameters (standard error), goodness-of-fit test and LRT statistic ($p$-value) based on accelerated $p$-max/$p$-min stable and single $p$-max/$p$-min stable models. Data sets are ground level ozone data and Stanford heart plant data.}
    \label{Realdata}
     \resizebox{\linewidth}{!}{ 
    \begin{tabular}{c|r|r|r|r} \mytoprule
    \multirow{2}{*}{-}& \multicolumn{2}{c|}{\hu{Ground Level Ozone}}& \multicolumn{2}{c}{Stanford Heart Plant} \\\cline{2-5}
     & accelerated $p$-max type & $p$-max stable distribution & accelerated $p$-min type & $p$-min stable distribution \\\hline
    $\widehat{\mu}_1$  & $2.74(0.34)$ & $4.52(0.01)$ & $7.26(0.15)$  & $5.31(0.19)$ \\\hline
    $\widehat{\sigma}_1$ & $1.75(0.27)$ & $0.15(<0.005)$  & $0.23(0.14)$  & $1.62(0.14)$ \\\hline
    $\widehat{\xi}_1$ & $-0.76(0.01)$  & $-0.24(0.02)$  & $0.30(0.33)$  & $-0.13(0.09)$\\\hline
    $\widehat{\mu}_2$ & $4.51(0.01)$ & -  & $5.38(0.26)$   & -\\\hline
    $\widehat{\sigma}_2$ & $0.14(<0.005)$ & - & $2.04(0.31)$   & -\\\hline
    $\widehat{\xi}_2$ & $-0.39(0.03)$ & - & $-0.30(0.10)$ & -\\\hline
    KS & $\mathbf{0.0316(0.69)}$ & $0.0484(0.18)$ & $\mathbf{0.0493(0.96)}$ & $0.0603(0.85)$ \\\hline
    CVM & $\mathbf{0.0777(0.71)}$ & $0.2827(0.15)$ & $\mathbf{0.0330(0.97)}$ & $0.0817(0.68)$ \\\hline
    AD & $\mathbf{0.4915(0.76)}$ & $1.8708(0.11)$ & $\mathbf{0.2167(0.99)}$ & $0.5570(0.69)$ \\\hline
    LRT & \multicolumn{2}{c|}{$25.7800(< 0.001)$} & \multicolumn{2}{c}{$12.9160(0.005)$} \\\mybottomrule
    \end{tabular}
    }
    \label{Table5}
\end{table}

\newpage
\appendix
\section{Appendix}\label{Appendix}
The proof of Theorem \ref{TheoremMaxMax} is given in Appendix \ref{AppendixA}, followed by Lemma \ref{lemmaA} \citep[Lemma 2.1]{Baracat2002} for the equivalent power normalization constants in Appendix \ref{AppendixB}. We end this section with Appendix \ref{Appendix C} for the proof of Theorem \ref{Theoremconsistency}. 

\subsection{Proof of Theorem \ref{TheoremMaxMax}} \label{AppendixA}
{Firstly, we rewrite the distribution of the power normalized maximum of $M_n$ defined in Eq.\eqref{MaxMax} as follows.} 
\begin{eqnarray} \label{MaxMax2}
&&\pk{\alpha_{2,n_2}\abs{M_n}^{ \beta_{2,n_2} } \sign(M_n) \le x} \notag \\ &=&\pk{ M_n \le \abs{ \frac{x}{\alpha_{2,n_2}} }^{ 1/\beta_{2,n_2} } \sign(x) } \notag \\
    &=& \pk{\alpha_{2,n_2}\abs{M_{1,n_1}}^{\beta_{2,n_2}} \sign(M_{1, n_1})\le x } \notag \\ &&\times \pk{ \alpha_{2,n_2}\abs{M_{2,n_2}}^{\beta_{2,n_2}} \sign(M_{2, n_2})\le x} \notag \\
    &=:& I_n \cdot I\!I_n. 
\end{eqnarray}
It follows by Eq.\eqref{ComponentPowerLimit} that
\begin{eqnarray} \label{II}
    I\!I_n \toweak H^{(2)}(x), \quad n_2\to\infty.
\end{eqnarray}
Now it remains to show the limit of $I_n$. 
   \begin{eqnarray*}
       I_n &=&\pk{M_{1,n_1} \le \abs{ \frac{x}{\alpha_{2,n_2}} }^{ 1/\beta_{2,n_2} } \sign(x)}  \\&=& \pk{ M_{1,n_1} \le \fracl{x_n}{\alpha_{1,n_1}}^{ 1/\beta_{1,n_1} } \sign(x)}, 
   \end{eqnarray*} 
   where $x_n =\alpha_{1,n_{1}} \fracl {|x|}{\alpha_{2,n_2}}^{\beta_{1,n_1} / \beta_{2,n_2}} =:\alpha_n\abs{x}^{\beta_n}$ with $\alpha_n, \beta_n$ given by Eq.\eqref{NormingMaxMax}. 
   
\underline{Case (i)}: Note that condition \eqref{ConditionNorming} holds for positive $A,B$. It follows further by uniform convergence in \eqref{ComponentPowerLimit} 
that
   $$I_n \toweak H^{(1)}(  A|x|^B \sign(x) )\quad \mbox{as}\ \min(n_1, n_2)\to\infty.$$
This, together with Eqs.\eqref{MaxMax2} and \eqref{II}, implies the claim. 

\underline{Case (ii)}: {Recalling Eqs.\eqref{MaxMax2} and \eqref{II}, it remains to show $I_n\to1$. Below, we will discuss the limit according to the supports of the six $p$-types distribution.} Indeed, the support of the distributions of the same $p$-types of $H_{1,\alpha}, H_{2,\alpha}, H_{5}$ is a subset of $(0, \infty)$, while the right-endpoints of $H_{3,\alpha}, H_6$, and $H_{4, \alpha}$ are zero and negative, respectively. 
\begin{itemize}
    \item[a).] For $H^{(2)}$  being one of the same $p$-types of $H_{1,\alpha}, H_{2,\alpha}, H_{5}$, while $H^{(1)}$  being one of the same $p$-types of $H_{3,\alpha}, H_{4,\alpha}, H_{6}$, we have for $x \ge 0$ 
    $$H^{(1)}( x_n \sign(x) ) =  H^{(1)}( \alpha_n |x|^{\beta_n} \sign(x) ) \ge H^{(1)}(0)= 1.$$
    The claim follows thus straightforwardly. 
    \item[b).] For $H^{(2)}$ being one of the same $p$-types of $H_{1,\alpha}, H_{2,\alpha}, H_{5}$, and $H^{(1)}$ being one of the same $p$-types of $H_{1,\alpha},  H_{5}$, we have for $A=\infty, 0\le B < \infty$ 
    $$\log x_n = \log \alpha_n + \beta_n \log x \to \infty\quad \mbox{as} \min(n_1, n_2)\to\infty $$
    holds for all $x>0$.     Noting that the limit in \eqref{ComponentPowerLimit} holds uniformly for $x\in\mathbb R$, the claim follows since $I_n \to 1$.     
    \item[c).] For $H^{(2)}$ being one of the same $p$-types of $H_{1,\alpha}, H_{2,\alpha}, H_{5}$, and $H^{(1)}$ being the same $p$-type of $H_{2,\alpha}$, we have for $A = \infty, 0 \le B < \infty$ or $A\ge \gamma(H^{(1)}) >0, B=0$ and any $x>0$
    $$\log x_n = \log \alpha_n + \beta_n \log x \to \log A\ge \log \gamma(H^{(1)}) \quad \mbox{as} \min(n_1, n_2)\to\infty. $$
    Similar argument as for $b)$ applies and  the claim follows since $I_n \approx H^{(1)}(x_n) \to 1$.
     \item[d).] For $H^{(1)}, H^{(2)}$ being one of the same $p$-types of $H_{3,\alpha}, H_{4,\alpha}, H_{6}$, we have for $0\le A\le -\gamma(H^{(1)}), \kai{B=0}$ or $A=0, 0\le B <\infty$ and any $x<0$
    $$x_n \sign(x) = -\alpha_n |x|^{\beta_n} \ge \gamma(H^{(1)}),$$
    which indicates that $I_n \approx H^{(1)}(x_n \sign(x)) \to 1$.
   \end{itemize}

\underline{Case (iii)}: 
Consider first that both $H^{(1)}$ and $ H^{(2)}$ are one of the same $p$-types of $H_{1,\alpha}, H_{2,\alpha}, H_{5}$ and thus let $x>0$ in the following. We rewrite $x_n$ as below
$$\log x_n= 
\beta_n \left(\frac{\log \alpha_n}{\beta_n} + \log x \right). 
$$ 
Since $\lim_{n \to \infty} (\log \alpha_n)/\beta_n=-C<0$, we have 
$x_n$ converges to $\infty$ for $x> x_0 =\e^C$, and $0$ for $x< x_0$. This, together with the left-continuity of distributions and uniform convergence, implies that 
$$I_n = H^{(1)}(x_n \sign(x)) = H^{(1)}(x_n) \to \I{x >x_0}.$$
Consequently, we obtain the claim for the first case. Similar argument applies for $H^{(1)}, H^{(2)}$ being one of the same $p$-types of $H_{3,\alpha}, H_{4,\alpha}, H_{6}$.
\\
\COM{Further, let $x_0=\e^C$ be in the interior of the support of $H^{(2)}$, we will consider following three cases:
\begin{enumerate}
    \item For $x>\e^C$, $\log(x_n) \to \infty$, then we have $x_n \sign(x) \to \infty$, which leads to $H^{(1)}(x_n \sign(x)) \to 1$.
    \item For $-\e^C<x<\e^C$, $\log(x_n) \to -\infty$, then we have $x_n \sign(x) \to 0$ which leads to $H^{(1)}(x_n \sign(x)) \to 0$.
    \item For $x<-\e^C$, $\log(x_n) \to \infty$, then we have $x_n \sign(x) \to -\infty$ which leads to $H^{(1)}(x_n \sign(x)) \to 0$.
\end{enumerate}
Hence we can conclude that for $x>\e^C$ we have $H^{(1)}(x_n\sign(x)) \to 1$. While for $x<\e^C$ will it will cause that $H^{(1)}(x_n\sign(x)) \to 0$, thus we have
 $$I_n \toweak \I{x > x_0}. $$
For $A=0$, $B=\infty$ and $H^{(1)}, H^{(2)}$ being one of the same $p$-types of $H_{3,\alpha}, H_{4,\alpha}, H_{6}$ under proper conditions of sample length, let $x_0=-\e^C$ be in the interior of the support of $H^{(2)}$. We can also obtain following three cases:
\begin{enumerate}
    \item For $x>\e^C$, $\log(x_n) \to \infty$, then we have $x_n \sign(x) \to \infty$, which leads to $H^{(1)}(x_n \sign(x)) \to 1$.
    \item For $-\e^C<x<\e^C$, $\log(x_n) \to -\infty$, then we have $x_n \sign(x) \to 0$ which leads to $H^{(1)}(x_n \sign(x)) \to 1$.
    \item For $x<-\e^C$, $\log(x_n) \to \infty$, then we have $x_n \sign(x) \to -\infty$ which leads to $H^{(1)}(x_n \sign(x)) \to 0$.
\end{enumerate}
Therefore, we can summarize that for $x>-\e^C$ we have $H^{(1)}(x_n\sign(x)) \to 1$. While for $x<-\e^C$ will it will cause that $H^{(1)}(x_n\sign(x)) \to 0$, thus we also have
$$I_n \toweak \I{x > x_0}. $$
\COM{
Firstly denote $\gamma_*(H^{(1)})=\inf\left\{x\in\mathbb R: H^{(1)}(x)>0\right\}$.
For $H^{(1)}, H^{(2)}$ being one of the same $p$-types of $H_{1,\alpha}, H_{2,\alpha}, H_{5}$, we have for $A=0, B = \infty$ and any $x>0$
    $$\log x_n= 
\beta_n(\frac{\log \alpha_n}{\beta_n} + \log x)= \beta_n(\log x - C) \to \left\{ 
\begin{array}{ll}
x_n>\gamma(H^{(1)}), &x>x_0,\\
x_n<\gamma_*(H^{(1)}), & x \le x_0
\end{array}
\right. $$
which will leads to
$$ \left\{ 
\begin{array}{ll}
\log(\gamma(H^{(1)})), &x>x_0,\\
\log(\gamma_*(H^{(1)})), & x \le x_0
\end{array}
\right.
$$
with $x_0$ being in the interior of the support of $H^{(2)}$. 
This will lead to 
$$I_n \to H^{(2)}(x) \I{x > x_0}. $$
The proof for $H^{(1)}, H^{(2)}$ being one of the same $p$-types of $H_{3,\alpha}, H_{4,\alpha}, H_{6}$ is similar with $x_0=-\expon{C}$, and thus the details are not repeated once more.
}
}
Consequently, we complete the proof of Theorem \ref{TheoremMaxMax}. 
\QED

\subsection{Equivalent Power Normalization Constants}\label{AppendixB}
Similar to \citet[Theorem 1.2.3 (Khintchine)]{leadbetter2012extremes}, the lemma below indicates the power limit type of invariance principle under power normalization.
\begin{lemma}[Lemma 2.1 in \citet{Baracat2002}]\label{lemmaA}
    Let $F_n$ be a sequence of distribution functions and $T_1$ be a non-degenerate distribution function. Suppose that $\alpha_n, \beta_n$ are positive  constant sequences  such that 
\begin{equation*}
    F_n\left( \left|\frac{x}{\alpha_n}\right|^{1/\beta_n} \sign(x)\right) \toweak T_1(x),
\end{equation*}
    then for some non-degenerate distribution function $T_2$ and positive constants $\alpha_n^\prime $ and $\beta_n^\prime $,
\begin{equation*}
     F_n\left( \left|\frac{x}{\alpha_n^\prime}\right|^{1/\beta_n^\prime} \sign(x)\right) \toweak T_2(x)
\end{equation*}
if and only if
\begin{equation*}
    \fracl{\alpha_n^{\beta_n^\prime}}{(\alpha_n^\prime)^{\beta_n}}^{1/\beta_n^\prime} \to A > 0 , \quad \frac{\beta_n}{\beta_n^\prime} \to B > 0
\end{equation*}
for some $A>0$ and $B>0$ and then
\begin{equation*}
    T_2(x)=T_1(Ax^B \sign(x)).
\end{equation*}
\end{lemma}

\subsection{Proof of Theorem \ref{Theoremconsistency}
}\label{Appendix C}
Before we show Theorem \ref{Theoremconsistency}, we recall Theorem 1 in \citet{SmithRichardL.1985Mlei} (stated below as Lemma \ref{Lemma:  Consistency}), which requires  Assumptions 1, 6-8 therein, stated below as Assumptions 1$\sim$4 for the density of the form in Eq.\eqref{Eq: Smith}. 
\begin{assumption}\label{assumption1}
    All second-order partial derivatives of $g(x;\bm{\phi})$ with respect to $x$ and $\bm{\phi}$ exist and are continuous in $0<x<\infty$ and $\bm{\phi}\in\Phi\subset \mathbb{R}^{p}$. Moreover, $c(\bm{\phi})=\alpha^{-1}\lim_{x \downarrow 0} g(x;\bm{\phi})$ exists, is positive and finite for each $\bm{\phi}$, and twice continuously differentiable as a function of $\bm{\phi}$.
\end{assumption}
\begin{assumption}\label{assumption2}
    Let $\EE_{\bm{\phi}}$ denote expectation with respect to $f(\cdot;0,\bm{\phi})$, suppose that for each $\bm{\phi}= (\phi_1,\ldots, \phi_p)$\\
 (a)
\begin{align*}
&\EE_{\bm{\phi}} \left \{ \left ( \frac{\partial }{\partial \phi_i}  \right ) \log f(X;0,\bm{\phi} ) \right \}=0,\quad i=1,\ldots, p\\
&\EE_{\bm{\phi}} \left \{ \left ( \frac{\partial }{\partial \phi_i} \right ) \log f(X;0,\bm{\phi} ) \left ( \frac{\partial }{\partial \phi_j}  \right ) \log f(X;0,\bm{\phi} ) \right \} \\
&\qquad =-\EE_{\bm{\phi}} \left \{\frac{\partial^2}{\partial \phi_i\partial \phi_j}\log f(X;0,\bm{\phi} )  \right \} =m_{ij}(\bm{\phi}),\quad i,j=1,\ldots, p;
\end{align*}

(b) for $\alpha>1$, 

\begin{align*}
&\EE_{\bm{\phi}} \left \{ \left ( \frac{\partial }{\partial x}  \right ) \log f(X;0,\bm{\phi} ) \right \}=0,\\
&-\EE_{\bm{\phi}} \left \{ \left ( \frac{\partial }{\partial x}  \right ) \log f(X;0,\bm{\phi} ) \left ( \frac{\partial }{\partial \phi_i}  \right ) \log f(X;0,\bm{\phi} ) \right \} \\
&\quad =\EE_{\bm{\phi}} \left \{ \frac{\partial^2}{\partial x\partial \phi_i} \log f(X;0,\bm{\phi} )  \right \}=m_{i0}(\bm{\phi})=m_{0i}(\bm{\phi});
\end{align*}

(c) further for $\alpha>2$,
\begin{align*}
\EE_{\bm{\phi}}\left [ \left \{ \left ( \frac{\partial }{\partial x}  \right )\log f(X;0,\bm{\phi} )  \right \}^{2}   \right ]&=-\EE_{\bm{\phi}} \left \{  \frac{\partial^2}{\partial x^2}\log f(X;0,\bm{\phi} )  \right \}\\&=m_{00}(\bm{\phi})>0.
\end{align*}
\end{assumption}

\begin{assumption}\label{assumption3}
    Let $h(x;\bm{\phi})$ be any of $\left (  \frac{\partial^2}{\partial x\partial \phi_i} \right ) \log g(x;\bm{\phi} ) $ or $\left (  \frac{\partial^2}{\partial \phi_i\partial \phi_j} \right ) \log g(x;\bm{\phi} ) $, then, as $\theta \to \theta_0$, $\bm{\phi} \to \bm{\phi}_0$,
    $$\EE_0\left | h(X-\theta;\bm{\phi}  ) - h(X-\theta_0;\bm{\phi}_0  )\right  | \to 0, $$ where $\EE_0$ is with respect to $f(\cdot;\theta_0,\phi_0 )$. If $\alpha>2$, we require the same of $h(x;\bm{\phi})=\left (  \frac{\partial^2}{\partial x^2} \right )  \log g\left ( x;\bm{\phi}  \right )$.
\end{assumption}

\begin{assumption}\label{assumption4}
    For each $\varepsilon>0$, $\delta>0$, there exists a function $h_{\varepsilon,\delta}$ such that 
    $$\left | \left (  \frac{\partial^2}{\partial x^2}  \right )\log g\left ( x;\bm{\phi}  \right )   \right | <\frac{\varepsilon }{x^{2} }+ h_{\epsilon,\delta}(y;\bm{\phi}_0  )$$ whenever $\lVert\bm\phi -\bm\phi_0\rVert<\delta$, $\left | x-y  \right |<\delta $, and $h_{\varepsilon,\delta}$ satisfies
    $$\int_{0}^{\infty } h_{\varepsilon,\delta}(x;\bm{\phi}_0  )  f(x;0,\bm{\phi}_0)dx <\infty.$$  
\end{assumption}

Given a random sample $X_1,\ldots,X_n$ from $X$ with density $f(\cdot;\theta,\bm{\phi})$ specified in Eq.\eqref{Eq: Smith}. Denote by $\ell_n(\theta,\bm\phi)$  the log-likelihood given by
\begin{equation*}
    \ell_n(\theta,\bm{\phi})=\sum_{i=1}^{n} \log f\left ( X_i;\theta ,\bm{\phi} \right ).
\end{equation*}
The pseudo-maximum likelihood estimator, when it exists, will be denoted by $(\widehat{\theta}_n,\widehat{\bm{\phi}}_n)$, the solution of the normal equations satisfying
\begin{equation*}
    \frac{\partial \ell_n \left (\theta,\bm{\phi}  \right ) }{\partial \theta}=0, \quad \frac{\partial \ell_n \left (\theta,\bm{\phi}  \right ) }{\partial \phi_i}=0,\quad i=1,\ldots,p.
\end{equation*}
Let $\bar{\bm{\phi}}_n$ denote the maximum likelihood estimator for $\bm{\phi}$ when $\theta=\theta_0$ is known, 
similarly, let $\bar{\theta}_n$ denote the maximum likelihood estimator for $\theta$ when $\bm{\phi}=\bm{\phi}_0$ is known. Define $\zeta_{n,\alpha}, \, n\ge 1$ to be 1 if $\alpha>1$ and $\log n$ if $\alpha=1$, and write $Y_n<_p r_n$ for random variables $\left\{Y_n\right\}$ and positive constants $\left\{r_n\right\}$ if 
\begin{equation*}
    \lim_{a \to \infty}\lim_{n \to \infty} \sup \mathbb{P}\left ( \left | Y_n \right | >ar_n \right )=0. 
\end{equation*}

\begin{lemma}[Theorem 1 in \citet{SmithRichardL.1985Mlei}]\label{Lemma:  Consistency}
    Suppose that Assumptions 1-4 above are satisfied for Eq.\eqref{Eq: Smith}. Let $\bm M$ be the strictly positive-definite matrix as $(m_{ij}(\bm{\phi}_0))_{i,j=0,\cdots,p}$ for $\alpha>2$, and $(m_{ij}(\bm{\phi}_0))_{i,j=1,\cdots,p}$ for $1<\alpha \leq 2$. Then there exist a sequence $(\widehat{\theta}_n,\widehat{\bm{\phi}}_n)$ of solutions to the likelihood equations such that
     \begin{equation*}
         \widehat{\theta}_n-\theta _{0} <_p (n\zeta_{n,\alpha } )^{-\frac{1}{2} },\quad  \widehat{\bm{\phi} }_n- \bm{\phi}_0<_p n^{-\frac{1}{2} }.
     \end{equation*}
     Moreover, if $\alpha=2$ we have
     \begin{equation*}
         \widehat{\theta }_n-\bar\theta _n <_p n^{-\frac{1}{2} }(\log n)^{-1}   ,\quad  \widehat{\bm{\phi} }_n- \bar{\bm{\phi}}_n<_p (n\log n)^{-\frac{1}{2} },
     \end{equation*}
     while if $1<\alpha<2$ we have
     \begin{equation*}
         \widehat{\theta }_n-\bar\theta _n <_p n^{-\frac{2}{\alpha }  +\frac{1}{2} }   ,\quad  \widehat{\bm{\phi} }_n- \bar{\bm{\phi}}_n<_p n^{-\frac{1}{\alpha } }.
     \end{equation*}
\end{lemma}

\COM{\begin{lemma}[Theorem 3 \Hu{(i) and (ii)} in \cite{SmithRichardL.1985Mlei}]\label{Lemma:  Asymptotic}
    Under the assumptions of Lemma \ref{Lemma:  Consistency} let $(\widehat{\theta}_n,\widehat{\bm{\phi}}_n)$ denote a sequence of maximum likelihood estimators satisfying the conclusions of Lemma \ref{Lemma:  Consistency}. \Hu{Here $c=c(\bm{\phi})$}.
    \begin{itemize}
        \item[(i).] If $\alpha>2$ then $\sqrt{n}(\widehat{\theta}_n-\theta_0,\widehat{\bm{\phi}}_n-\bm{\phi}_0)$ converges in distribution with mean $\bm0$ and covariance matrix $\bm{M}^{-1}$, where $\bm{M}$ is as in Lemma \ref{Lemma:  Consistency}.
        \item[(ii).] If $\alpha=2$ then $\left \{ \sqrt{nc\log n}(\widehat{\theta}_n-\theta_0),\sqrt{n}(\widehat{\bm\phi}_n-\bm{\phi}_0) \right \} $ converges in distribution to a normal random vector with covariance matrix of the form 
        \begin{equation*}
            \begin{bmatrix}
             1&0 \\
             0&\bm{M}^{-1} 
            \end{bmatrix},
        \end{equation*}
        where $\bm{M}$ is as in Lemma \ref{Lemma:  Consistency}.
    \end{itemize}
\end{lemma}
}
\textit{\textbf{Proof of Theorem \ref{Theoremconsistency}.}}
We will show that the density given by Eq.\eqref{acc-l-min-pdf} satisfies all the four assumptions in Lemma \ref{Lemma:  Consistency}. The claim thus follows straightforwardly. \\
\underline{Verification of Assumption 1}. Clearly,  all second-order partial derivatives of $g(x;\bm{\phi})$ defined in Eq.\eqref{g-acc-l-min} with respect to $x$ and $\bm{\phi} = (\sigma_1,\alpha_1,\sigma_2,\alpha_2)\in \Phi$ exist and are continuous in $0<x<\infty$ and $\bm{\phi}\in\Phi$. In addition,
\begin{eqnarray*}
    c(\bm{\phi} ) =\alpha_1^{-1}\lim_{x \to 0}g(x;\sigma_1,\alpha_1,\sigma_2,\alpha_2)=\left(\frac{1}{\sigma _1\alpha _1}
    \right)^{\alpha_1}>0. 
\end{eqnarray*}
 Therefore, {Assumption 1} follows. \\
\underline{Verification of Assumption 2}.
We consider the first-order derivative of the pdf $\underline h$ as below. 
\begin{eqnarray*}
  \frac{\partial \underline{h}(x ;0,\bm{\phi})}{\partial \sigma_1} &=&   \frac{\alpha_1}{\sigma_1}\fracl{x}{\sigma_1\alpha_1}^{\alpha_1}\left [ \frac{1}{\sigma _1} \left ( \frac{x}{\sigma_1\alpha_1}  \right )^{\alpha_1-1}+ \frac{1}{\sigma _2} \left ( \frac{x}{\sigma_2\alpha_2}  \right )^{\alpha_2-1} -\frac{\alpha_1}x \right ] \\
  &&\times \exp\left[-\left ( \frac{x}{\sigma _1\alpha _1}  \right ) ^{\alpha_1}-\left ( \frac{x}{\sigma _2\alpha _2}  \right ) ^{\alpha_2}\right].    
\end{eqnarray*}
Note that there exists a generate constant $K$ (which varies from line to line) such that
\begin{eqnarray*}
    \left | \frac{\partial \underline{h}(x ;0,\bm{\phi})}{\partial \sigma_1} \right |&\le& K\left(x^{2\alpha_{1}-1}+x^{\alpha_{1}+\alpha_{2}-1}+x^{\alpha_{1}-1}\right) \exp\left[-\left ( \frac{x}{\sigma _1\alpha _1}  \right ) ^{\alpha_1}-\left ( \frac{x}{\sigma _2\alpha _2}  \right ) ^{\alpha_2}\right]
\end{eqnarray*}
which is integrable on $(0,\infty)$ for $1<\alpha_1 \le \alpha_2$. Thus, we have
\begin{equation*}
    \int_{0}^{\infty } \frac{\partial \underline{h}(x ;0,\bm{\phi})}{\partial \sigma _1}dx= \frac{\partial }{\partial \sigma_1} \int_{0}^{\infty}\underline h(x;0,\bm{\phi})dx=0
\end{equation*}
by dominated convergence theorem.
Similar argument applies  for $\sigma_2$ since $\sigma_1$ and $\sigma_2$ are exchangeable.
\\
Next, we consider the first-order parital derivative of $\underline h$ with respect to $\alpha_1$.  Note that 
\begin{eqnarray*}
\frac{\partial \underline{h}(x ;0,\bm\phi)}{\partial \alpha_1} &=& \fracl{x}{\sigma_1\alpha _1}^{\alpha_1} \Bigg\{ \left[ \log\left( \frac{\sigma _1\alpha_1}{x} \right) + 1 \right] \left[\frac{1}{\sigma _1} \left( \frac{x}{\alpha_1\sigma_1} \right)^{\alpha_1-1} + \frac{1}{\sigma _2} \left( \frac{x}{\alpha_2\sigma_2} \right)^{\alpha_2-1}\right] \\
&&+ \frac{\alpha _1}{x}\left[\log\left( \frac{x}{\sigma _1\alpha_1} \right) + \frac{1}{\alpha _1} - 1 \right] \Bigg\}   \exp\left[-\left( \frac{x }{\sigma _1\alpha _1} \right)^{\alpha_1} - \left( \frac{x }{\sigma _2\alpha _2} \right)^{\alpha_2}\right],
\end{eqnarray*}
which is integrable on $(0,\infty)$ (recalling Gumbel distribution and Gamma distribution exists with finite mean). Therefore, using dominated convergence theorem, we have
\begin{equation*}
    \int_{0}^{\infty } \frac{\partial \underline{h}(x ;0,\bm\phi)}{\partial \alpha_1}dx= \frac{\partial }{\partial \alpha_1} \int_{0}^{\infty}\underline h(x;0,\bm\phi)dx=0.
\end{equation*}
Similarly, one may show all second-order derivatives $\partial^2/{\partial \phi_i\partial \phi_j} \underline{h}(x ;0,\phi)$
for $(\phi_i,\phi_j)$ being one of $(\sigma_1,\sigma_2), (\sigma_1,\alpha_1),\\(\sigma_1,\alpha_2), (\sigma_2,\alpha_1),(\sigma_2,\alpha_2)$ exist with integrable bounds and thus using dominated convergence theorem, Assumption 2 holds with variance-covariance matrix $\bm M$ positive definite. 
\\
\underline{Verification of Assumption 3}. 
We have 
\begin{eqnarray*}
\log g(x;\bm{\phi}) &=& \log\left [ \frac{1}{\sigma _1} \left ( \frac{1}{\sigma_1\alpha_1}  \right )^{\alpha_1-1}+ \frac{1}{\sigma _2} \left ( \frac{1}{\sigma_2\alpha_2}  \right )^{\alpha_2-1}  x^{\alpha_2-\alpha_1} \right ] \\
&&- \left ( \frac{x}{\sigma _1\alpha _1}  \right ) ^{\alpha_1}-\left ( \frac{x}{\sigma _2\alpha _2}  \right ) ^{\alpha_2},\\
    \frac\partial{\partial x}\log g(x;\bm{\phi}) &=&  \frac{(\alpha_2-\alpha_1) x^{\alpha_2-\alpha_1 -1}}{x^{\alpha_2 - \alpha_1}+  \frac{\sigma_2}{\sigma _1}  \frac{(\sigma_2\alpha_2)^{\alpha_2-1}}{(\sigma_1\alpha_1)^{\alpha_1-1}}} - \frac{\alpha_1 x^{\alpha_1 -1}}{(\sigma_1\alpha_1)^{\alpha_1}} - \frac{\alpha_2 x^{\alpha_2 -1}}{(\sigma_2\alpha_2)^{\alpha_2}}, \\
        \frac{\partial^2}{\partial x^2}\log g(x;\bm{\phi}) &=&   - 
        \left(\frac{(\alpha_2-\alpha_1) x^{\alpha_2-\alpha_1 -1}}{x^{\alpha_2 - \alpha_1}+  \frac{\sigma_2}{\sigma _1}  \frac{(\sigma_2\alpha_2)^{\alpha_2-1}}{(\sigma_1\alpha_1)^{\alpha_1-1}}}\right)^2 + \frac{(\alpha_2-\alpha_1) (\alpha_2-\alpha_1-1)x^{\alpha_2-\alpha_1 -2}}{x^{\alpha_2 - \alpha_1}+  \frac{\sigma_2}{\sigma _1}  \frac{(\sigma_2\alpha_2)^{\alpha_2-1}}{(\sigma_1\alpha_1)^{\alpha_1-1}}} \\
        &&
        - \frac{\alpha_1(\alpha_1 -1) x^{\alpha_1 -2}}{(\sigma_1\alpha_1)^{\alpha_1}} - \frac{\alpha_2 (\alpha_2-1)  x^{\alpha_2 -2}}{(\sigma_2\alpha_2)^{\alpha_2}}
\end{eqnarray*}
are bounded as $x\downarrow0$ respectively for $0< \alpha =\alpha_1\le \alpha_2; 1<\alpha =\alpha_1 = \alpha_2$ or $1<\alpha =\alpha_1 <\alpha_2 -1$; $2<\alpha =\alpha_1 = \alpha_2$ or $2< \alpha =\alpha_1 <\alpha_2 -2$. Thus, Assumption 3 follows.
\\
\underline{Verification of Assumption 4}. The assumption holds for the case with $\alpha_1 = \alpha_2$, following the same argument for the Weibull distribution. Below, we consider only that $2<\alpha_1 < \alpha_2$. First of all, we have with $\alpha_m = \min\{\alpha_2 - \alpha_1 - 2, \alpha_1 -2\}, \alpha_M = \max\{2(\alpha_2 - \alpha_1 - 1), \alpha_2 -2\}$
\begin{eqnarray*}
    \abs{\frac{\partial^2}{\partial x^2}\log g(x;\bm{\phi})} &\le& K\left[x^{2(\alpha_2 - \alpha_1 -1)} + x^{(\alpha_2 - \alpha_1 - 2)} + x^{\alpha_1 -2} + x^{\alpha_2 -2}\right]\\
    &\le & 
    \left\{\begin{array}{ll}
        K x^{\alpha_{m}}, & 0\le x\le 1, \\
        K x^{\alpha_{M}}, &x>1
        \end{array}
        \right.
\end{eqnarray*}
hold for all $\phi_0 =(\alpha_{10},\sigma_{10}, \alpha_{20}, \sigma_{20})$ such that $\alpha_m < \alpha_{m0} <\alpha_{M0} < \alpha_M$ and a generate constant $K$.\\
For $\bm\phi$ such that $\lVert\bm\phi -\bm\phi_0\rVert < \delta$ with a given $\delta>0$, we take $\overline\alpha_0 = (\alpha_{m0}+\alpha_{M0})/2$ and $x_\varepsilon=\left(\varepsilon/K\right)^{1/(\overline\alpha_0+2)}$  which is assumed to be less than or equal to 1. We have 
\begin{eqnarray*}
    \abs{\frac{\partial^2}{\partial x^2}\log g(x;\bm{\phi})}\le\frac{\varepsilon}{x^2}+{h}'_{\varepsilon,\delta}(x;\bm\phi_0)
\end{eqnarray*}
where
\begin{eqnarray*}
    {h}'_{\varepsilon,\delta}(x;\bm\phi_0)={\begin{cases}
 0,&0<x<x_\varepsilon, \\
 K x^{\alpha_{m}} ,&x_\varepsilon\le x\le 1,\\
 K x^{\alpha_{M}} ,&x >1.
\end{cases}}
\end{eqnarray*}
This function is bounded as $x \to 0$. If we define $h_{\varepsilon,\delta}(y;\bm\phi_0)$ to be
\begin{eqnarray*}
    \sup\left \{ {h}'_{\varepsilon,\delta}(x;\bm\phi_0):\left | x-y \right |<\delta   \right \},
\end{eqnarray*}
Assumption 4 is satisfied.

Therefore, a straightforward application of Lemma \ref{Lemma: Consistency} completes the proof of Theorem \ref{Theoremconsistency}, i.e., we obtain the consistency of pseudo-maximum likelihood estimation. 
\QED

\bibliographystyle{apalike}

\end{document}